\documentclass[11pt]{article}

\usepackage{amsmath,amsthm,amsfonts,amssymb,amscd}
\usepackage{subfigure}
\usepackage{graphicx,bm,color}
\usepackage{algorithm}
\usepackage[algo2e]{algorithm2e}
\usepackage{stmaryrd}
\usepackage{lineno}
\usepackage{tikz}
\usepackage[colorlinks,allcolors=blue]{hyperref}

\definecolor{gray}{gray}{0.6}

\theoremstyle{plain}

\graphicspath{{./}}

\def\n{{\bm n}}

\begin{document}

\title{Multiscale Methods for Model Order Reduction of Non Linear Multiphase Flow Problems}

\author{Gurpreet Singh \and Wingtat Leung \and Mary F. Wheeler}

\date{\today}
\maketitle


\begin{abstract}
Numerical simulations for flow and transport in subsurface porous media often prove computationally prohibitive due to property data availability at multiple spatial scales that can vary by orders of magnitude. A number of model order reduction approaches are available in the existing literature that alleviate this issue by approximating the solution at a coarse scale. We attempt to present a comparison between two such model order reduction techniques, namely: (1) adaptive numerical homogenization and (2) generalized multiscale basis functions. We rely upon a non-linear, multi-phase, black-oil model formulation, commonly encountered in the oil and gas industry, as the basis for comparing the aforementioned two approaches. An expanded mixed finite element formulation is used to separate the spatial scales between non-linear, flow and transport problems. To the author's knowledge this is the first time these approaches have been described for a practical non-linear, multiphase flow problem of interest. A numerical benchmark is setup using fine scale property information from the 10$^{th}$ SPE comparative project dataset for the purpose of comparing accuracies of these two schemes. An adaptive criterion is employed in by both the schemes for local enrichment that allows us to preserve solution accuracy compared to the fine scale benchmark problem. The numerical results indicate that both schemes are able to adequately capture the fine scale features of the model problem at hand.
\end{abstract}





\newcommand{\bs}[1]{\boldsymbol{#1}}

\section{Introduction}
Since the advent of numerical reservoir simulations, oil-field operations have received a substantial boost in confidently predicting recovery estimates and determining operational choices during the deployment of a specific recovery technology. During the screening stage, numerical reservoir models are built and simulations are run to determine feasibility of a number of injection/production scenarios. This requires that the reservoir simulation to be both accurate and time efficient. Furthermore, for an oil and gas reservoir already in production, it is necessary to determine reservoir parameters (permeability, porosity, etc.) with increasing certainty. As the field matures, more efficient practices such as production stimulation, chemical treatment, or infill drilling can substantially increase the life of the reservoir and consequently the overall recovery.

Uncertainty quantification (UQ) and parameter estimation have been used in this respect to predict reservoir parameters incorporating well logs, production data and geological models during history matching. These frameworks either rely upon proxy models which approximate the flow physics in the reservoir, or a field scale reservoir simulator as the forward model driving engine for generating multiple realizations. The proxy models often suffer from inadequate representation of flow physics for example, a single phase flow model is often used as a proxy for more involved fluid description such as gas flooding that necessitates the use of an equation of state (EOS) compositional flow model. However, at the other end using a full-field reservoir simulator is prohibitively expensive since the number of realizations required by UQ for parameter estimation from a reservoir simulator, are of the order of 100  \cite{WRCR:WRCR20632}. Additionally, the field data is available at different spatial scales due to differences in observations techniques.  For well-logs this spatial scale is a few feet whereas data inferred from seismic recordings is of the order of 100 feet.

Highly heterogeneous medias often occur in oil reservoirs. The flow simulation will be strongly affected by these multiscale features of the reservoir. Features with different length scales will affect the flow solution differently. For example, fractures with length much larger than the coarse grid size will affect the solution globally. Periodic oscillating media with period much smaller than the coarse grid size or fractures with length much smaller than the coarse grid size will gives a local effect on the flow simulation.  There are a variety of techniques for handling multiscale features such as upscaling method and multiscale finite element methods. Upscaling methods,  \cite{durlofsky1991numerical,wu2002analysis,bourgeat1984homogenized,amaziane2010homogenization,chen2003coupled}, usually upscale the media properties and use properties to solve the flow problem on the coarse- grid. Multisacle methods, \cite{aarnes2004use,hou1997multiscale,jenny2003multi,arbogast2007multiscale,chen2003mixed}, construct coarse-grid basis functions and solve the problem by using these multiscale basis functions. Two major approaches are considered in this paper for upscaling, one based on two-scale homogenization theory and a second on generalized mulitscale finite elements.  Both of these schemes are described below.  Our objective in this paper is to compare the application of these two schemes to a realistic black oil model.

Two-scale homogenization theory is a mathematically consistent, theoretical framework that has been used by several others \cite{allaire1992homogenization,mikelic2006rigorous,jikov2012homogenization,bensoussan1978asymptotic, amaziane2006effective} to derive effective equations for a variety of problems. For a reservoir domain characterized by a length scale $L$, two-scale homogenization theory makes basic assumptions: (1) existence of an identifiable Representative Elementary Volume (REV) (or a period) with a characteristic length scale $l$, and (2) a scale separation between these aforementioned two length scale or $\varepsilon = l/L<<1$. The first assumption effectively means that a parameter or reservoir property that oscillates at the fine spatial scale ceases to do so at some coarse spatial scale. For further details regarding homogenization theory; both rigorous and formal, the reader is referred to \cite{yerlan} and the citations therein.

In adaptive homogenization, a local numerical homogenization is performed to reduce computational costs. More precisely a formal two-scale homogenization framework is utilized to derive effective equations at the coarse scale for a given fine scale flow model formulation. Adaptive mesh refinement (AMR) is then employed using an Enhanced Velocity Mixed Finite Element Method (EV MFE) \cite{wheeler2002enhanced,thomas2011enhanced} to recover fine scale flow features at the saturation front. Our approach shares similarities with the adaptive \cite{chung2016adaptive} and hybrid \cite{aarnes2006adaptive} multiscale methods in that a local enrichment is performed at the saturation front. These latter approaches rely upon enriching local multiscale basis functions as opposed to our approach of local mesh refinement in the region of interest. The saturation dependent mobilities are treated using an expanded mixed finite element formulation as discussed in Section \ref{sec:expanded}.

Generalized Multiscale Finite Element Method (GMsFEM), introduced in \cite{efendiev2013generalized}, uses local spectral problems to construct additional local basis functions to capture the effect
of these global features. In a GMsFEM, the computational process is divided into offline and online stages. A snapshot space and an offline multiscale space are constructed in the offline stage. The snapshot
space contains the functions which capture the multiscale properties of the solution. The offline multiscale space is a low dimension space of the snapshot space which gives a good approximate of the solution.
In the online stage, the multiscale space are used to compute the numerical solution with given boundary conditions and source functions. Sometimes, in the online stage, one may need to modify or select some of the basis functions in the offline space to construct an online multiscale space. We remark that the offline procedures are only preformed once for an fixed media. The offline multiscale basis function can
be used for the problems with different boundary conditions and different source functions.

In this work, a generalized mixed multiscale method is introduced for the black oil model problem, as an extension of the framework in \cite{chung2015mixed}. In \cite{chung2015mixed}, incompressible two phase flow was considered and and the multiscale basis functions were constructed to approximate the total velocity.  In the black oil model, the fluid compressibility effects total mobility which can lead to large changes in total velocity with time.  This requires that the basis functions for approximating total velocity be updated frequently.  To avoid this difficulty, an expanded mixed formulation is employed and the basis functions are constructed to approximate a pseudo velocity which is the product of the absolute permeability and the gradient of the reference pressure. Since the absolute permeability doel not vary in time, there is no need to update the multiscale basis functions during simulation runtime.

The outline of the paper is as follows:  in Section \ref{sec:formulation} we describe a three-phase black oil model.  In Section \ref{sec:expanded} we present the expanded mixed variational formulation used in the formulation of the aforementioned two model order reduction or upscaling approaches. In Sections \ref{sec:adaptive} and \ref{sec:GMMM}  adaptivie numerical homoenization and generalized mixed multiscale methods are introduced respectively.  Computational results are discussed in Section \ref{sec:results} followed by a summary in Section \ref{sec:summary}.

\section{Model Formulation} \label{sec:formulation}

In this section, we describe a three-phase, black oil model formulation with a local equilibrium model to represent live oil description in the oil and gas industry. We consider three component oil, water, and gas and label them with indices 1, 2, and 3 respectively in three fluid phases namely oil (o), water (w), and gas (g). The conservation equation for a component $i$ in phase $\alpha$ can then be written as,

\begin{equation}
\frac{\partial}{\partial t} \left(\phi \rho_{\alpha} S_{\alpha} x_{i\alpha} \right) + \nabla \cdot x_{i\alpha} \rho_{\alpha} u_{\alpha} = q_{i\alpha} + R_{i\alpha}, \quad i=1, 2, 3 \text{, and } \alpha = o,w,g.
\label{eqn:cons_comp}
\end{equation}

Here, $\phi$ is the fluid fraction (or porosity), and $\rho_{alpha}$, $S_{\alpha}$, and $\rho_{\alpha}$ the density, saturation and Darcy velocity of phase $\alpha$, respectively. Further, $x_{i\alpha}$ is the concentration of component $i$ in phase $\alpha$, $q_{i\alpha}$ the source/sink term, and $R_{i\alpha}$ is the mass exchange of component $i$ between phases due to local equilibrium description. The Darcy velocity of a phase $\alpha$ is given by,

\begin{equation}
u_{\alpha} = -K \frac{k_{r\alpha}}{\mu_{\alpha}} \nabla p_{\alpha}
\label{eqn:darcy}
\end{equation}
Here, $K$ is the absolute permeability, $k_{r\alpha}$ is the relative permeability of phase $\alpha$ defined as a function of phase saturations $S_{\alpha}$, and $\mu_{\alpha}$ and $p_{\alpha}$  the viscosity and pressure o of phase $\alpha$, respectively. In what follows, we define phase mobility as,
\begin{equation}
\lambda_{\alpha} = \frac{k_{r\alpha}}{\mu_{\alpha}},
\end{equation}
to reduce the notation. The black oil model assumes, that the gas component can exist only in the oil and gas phases whereas the oil and water component only exist in the oil and water phases, respectively. This results in the following constraint equations,
\begin{equation}
\begin{aligned}
\sum_{\alpha} R_{i\alpha} = 0\\
x_{1o} + x_{3o} = 1, \quad x_{2o} = 0,\\
x_{1w} = 0, \quad x_{2w} = 1, \quad x_{3w} = 0,\\
x_{1g} = 0, \quad x_{2g} = 0, \quad x_{3g} = 1.
\end{aligned}
\label{eqn:compcon}
\end{equation}
Please note that $x_{i\alpha}$ can be defined either as a mass fraction or mole fraction in the above description by consistently choosing phase densities on a mass or molar basis per unit phase volume. Additionally, the saturation and capillary pressure constraints are given by,

\begin{equation}
\begin{aligned}
\sum_{\alpha}  S_{\alpha} = 1,\\
P_{cow}(S_{w}) = P_{o} - P_{w},\\
P_{cgo}(S_{o}) = P_{g} - P_{o}.
\end{aligned}
\label{eqn:satcapcon}
\end{equation}
Here, $P_{cow}$ and $P_{cgo}$ are the capillary pressures corresponding to the oil-water and gas-oil interactions. Summing Eqn. \eqref{eqn:cons_comp} over all the phases we obtain the overall component conservation equations for the oil, water, and gas components as follows,

\begin{equation}
\begin{aligned}
\frac{\partial}{\partial t} \left( \phi \rho_{o} S_{o} x_{1o} \right) + \nabla \cdot \left(\rho_{o} x_{1o} u_{o} \right) = q_{1}\\
\frac{\partial}{\partial t} \left(\phi \rho_{w} S_{w} \right) + \nabla \cdot \left(\rho_{w} u_{w} \right) = q_{2}\\
\frac{\partial}{\partial t} \left(\phi \rho_{o} S_{o} x_{3o} +  \phi \rho_{g} S_{g} \right) + \nabla \cdot \left(\rho_{g} u_{g} + \rho_{o} x_{3o} u_{o} \right) = q_{3}
\end{aligned}
\end{equation}

This results in the familiar black oil formulation, with the solution gas defined by $R_{g} = x_{3o}$. For convenience we define the total component concentrations to replace the phase saturations as,
\begin{equation}
\begin{aligned}
N_{o} &= \rho_{o}S_{o}x_{1o}\\
N_{w} &= \rho_{w} S_{w}\\
N_{g} &= \left(\rho_{o} S_{o} x_{3o} + \rho_{g}S_{g}\right)
\end{aligned}
\label{eqn:consatrel}
\end{equation}
The phase saturations can then be evaluated using,
\begin{equation}
\begin{aligned}
S_{o} &= \frac{N_{o}}{\rho_{o} (1-R_{g})},\\
S_{w} &= \frac{N_{w}}{\rho_{w}},\\
S_{g} &= \frac{1}{\rho_{g}}\left(N_{g} - N_{o}\frac{R_{g}}{1-R_{g}}\right).
\end{aligned}
\label{eqn:satconrel}
\end{equation}
The saturation constraint is then expressed as,
\begin{equation}
\frac{N_{o}}{\rho_{o} (1-R_{g})} + \frac{N_{w}}{\rho_{w}} + \frac{1}{\rho_{g}}\left(N_{g} - N_{o}\frac{R_{g}}{1-R_{g}}\right) = 1
\label{eqn:satcon}
\end{equation}

%
%
%
%
%
%
%

\subsection{Black Oil Formulation}
With all the preliminaries defined, we now write the system of partial differential equations associated with the black oil formulation on the domain $\Omega \times (0,T]$ as,
\begin{equation}
\begin{aligned}
\frac{\partial}{\partial t} \left(\phi N_{o} \right) + \nabla \cdot \left(\rho_{o} x_{1o} u_{o} \right) = q_{1},\\
\frac{\partial}{\partial t} \left(\phi N_{w} \right) + \nabla \cdot \left(\rho_{w} u_{w} \right) = q_{2},\\
\frac{\partial}{\partial t} \left(\phi N_{g} \right) + \nabla \cdot \left(\rho_{g} u_{g} + \rho_{o} R_{g} u_{o} \right) = q_{3},
\end{aligned}
\label{eqn:system}
\end{equation}
with the phase Darcy velocity,
\begin{equation}
u_{\alpha} = -K \frac{k_{r\alpha}}{\mu_{\alpha}} \nabla p_{\alpha}.
\label{eqn:darcy2}
\end{equation}
The relative permeability ($k_{r\alpha}$) and solution gas ($R_{g}$) are defined as monotonic functions of phase saturations ($S_{\alpha}$) and reference phase pressure ($P_{ref}$), respectively. Further the constraint equations are,
\begin{equation}
\begin{aligned}
\frac{N_{o}}{\rho_{o} (1-R_{g})} + \frac{N_{w}}{\rho_{w}} + \frac{1}{\rho_{g}}\left(N_{g} - N_{o}\frac{R_{g}}{1-R_{g}}\right)& = 1\\
P_{cow}(S_{w})& = P_{ref} - P_{w}\\
P_{cgo}(S_{o}) &= P_{g} - P_{ref}
\end{aligned}
\label{eqn:constraint}
\end{equation}

Equations \eqref{eqn:system} thru \eqref{eqn:constraint} gives us a well determined system in the pressure ($P_{ref}$, $P_{w}$, $P_{g}$), concentration ($N_{o}$, $N_{w}$, $N_{g}$), Darcy velocity ($u_{\alpha}$, $\alpha=o,w,g$), and local equilibrium ($\nu$, $x_{1o}$) unknowns. Further, the boundary and initial conditions are given by,
\begin{eqnarray}
u_{\alpha} \cdot \n = 0 \text{~on~} \partial \Omega \times (0,T]\\
P_{ref} = P_{ref}^{0}, \quad N_{\alpha} = N_{\alpha}^{0},  \text{~at~} \Omega \times \{t=0\}
\end{eqnarray}
Please note that a no flow boundary condition is considered for convenience in model formulation and more general boundary conditions can also be treated.

\section{Expanded Mixed Variational Formulation}\label{sec:expanded}
In this section, we briefly present two, semi-discrete (space discrete, time continuous), expanded mixed variational forms of the black oil model formulation for adaptive numerical homogenization and generalized multiscale basis functions approaches.

\subsection{Expanded Mixed Form for Adaptive Numerical Homogenization}
We define velocity and pressure (or concentration) spaces as $\bs{V} = \{\bs{v} \text{~in~} H(\mathrm{div};\Omega) : \bs{v}\cdot \n = 0 \text{~on~} \partial\Omega\}$ and $W \equiv L^{2}(\Omega)$, respectively. Further, we define component fluxes for the oil, water, and gas components as,
\begin{equation}
\begin{aligned}
\bs{F}_{o}& = \rho_{o} x_{1o} u_{o},\\
\bs{F}_{w}& = \rho_{w} u_{w},\\
\bs{F}_{g}& = \rho_{g} u_{g} + \rho_{o}R_{g}u_{o}.
\end{aligned}
\end{equation}
The expanded mixed variational problem is: Find $N_{\alpha} \in W$, $P_{o}\in W$, $\bs{F}_{\alpha}\in V$, and $\tilde{\bs{F}}_{\alpha}\in V$ such that,

\begin{equation}
\begin{aligned}
\left(\frac{\partial}{\partial t} \left(\phi N_{o} \right),w\right)_{\Omega} + \left(\nabla \cdot \bs{F}_{o},w\right)_{\Omega}= (q_{1},w)_{\Omega},\\
\left(\frac{\partial}{\partial t} \left(\phi N_{w}\right),w \right)_{\Omega} + \left(\nabla \cdot \bs{F}_{w},w\right)_{\Omega} = (q_{2},w)_{\Omega},\\
\left(\frac{\partial}{\partial t} \left(\phi N_{g}\right),w \right)_{\Omega} + \left(\nabla \cdot \bs{F}_{g},w\right)_{\Omega} = (q_{3},w)_{\Omega},
\end{aligned}
\end{equation}
and
\begin{equation}
\begin{aligned}
\left(K^{-1}\tilde{\bs{F}}_{o},\bs{v}\right)_{\Omega} - \left(P_{o},\nabla \cdot \bs{v}\right)_{\Omega} &= 0,\\
\left(K^{-1}\tilde{\bs{F}}_{w},\bs{v}\right)_{\Omega} - \left(P_{o},\nabla \cdot \bs{v}\right)_{\Omega} &= -\left(P_{cow},\nabla \cdot \bs{v}\right)_{\Omega},\\
\left(K^{-1}\tilde{\bs{F}}^{1}_{g},\bs{v}\right)_{\Omega} - \left(P_{o},\nabla \cdot \bs{v}\right)_{\Omega} &= 0,\\
\left(K^{-1}\tilde{\bs{F}}^{2}_{g},\bs{v}\right)_{\Omega} & = \left(P_{cgo},\nabla \cdot \bs{v}\right)_{\Omega} ,
\end{aligned}
\label{eq:expanded_flux}
\end{equation}
and
\begin{equation}
\begin{aligned}
\left(\bs{F}_{o},\bs{v}\right)_{\Omega} &= \left(\lambda_{o} \tilde{\bs{F}}_{o},v\right)_{\Omega},\\
\left(\bs{F}_{w},\bs{v}\right)_{\Omega} &= \left(\lambda_{w} \tilde{\bs{F}}_{w},v\right)_{\Omega},\\
\left(\bs{F}_{g},\bs{v}\right)_{\Omega} &= \left(\lambda^{1}_{g} \tilde{\bs{F}}^{1}_{g},v\right)_{\Omega} + \left(\lambda^{2}_{g} \tilde{\bs{F}}^{2}_{g},v\right)_{\Omega}.
\label{eq:numerical_flux}
\end{aligned}
\end{equation}
Here, $w\in W$ and $\bs{v} \in \bs{V}$ and the mobilities $\lambda$ are defined as,
\begin{equation}
\begin{aligned}
\lambda_{o} &= \rho_{o} x_{1o} \frac{k_{ro}}{\mu_{o}},\\
\lambda_{w} &= \rho_{w} \frac{k_{rw}}{\mu_{w}},\\
\lambda^{1}_{g} & = \left(\rho_{g}\frac{k_{rg}}{\mu_{g}} + \rho_{o} R_{g} \frac{k_{ro}}{\mu_{o}}\right),\\
\lambda^{2}_{g} & = \rho_{g}\frac{k_{rg}}{\mu_{g}}.
\end{aligned}
\end{equation}
$\tilde{\bs{F}}_{\alpha}$ is the auxiliary velocity commonly used in the expanded mixed formulations to avoid inversion of zero corresponding to the phase relative permeability at irreducible (or residual) saturations.

\subsection{Expanded Mixed Form for Generalized Multiscale Basis Functions}
In this subsection, we will introduce the expanded mixed form for the generalized multiscale method. To solve a solution for the black oil model, we just need to seek for two of the saturations, one of the pressure and one of the velocities. In the generalized multiscale method, we are seeking for the water saturation, the gas saturation, the oil pressure and oil pseudo velocity. That is, we will find $(S_w(\cdot,t),S_g(\cdot,t),P_o(\cdot,t),\bs{\tilde{F}}_o(\cdot,t))\in W_{H,h}\times W_{H,h}\times W_{H,h}\times V_{H,h}$ such that
\begin{equation}
\begin{aligned}
\left(\frac{\partial}{\partial t} \left(\phi N_{o}(S_w,S_g,P_o) \right),w\right)_{\Omega} + \left(\nabla \cdot \bs{F}_{o},w\right)_{\Omega}&= (q_{1},w)_{\Omega},\\
\left(\frac{\partial}{\partial t} \left(\phi N_{w}(S_w,S_g,P_o)\right),w \right)_{\Omega} + \left(\nabla \cdot \bs{F}_{w},w\right)_{\Omega}&= (q_{2},w)_{\Omega},\\
\left(\frac{\partial}{\partial t} \left(\phi N_{g}(S_w,S_g,P_o)\right),w \right)_{\Omega} + \left(\nabla \cdot \bs{F}_{g},w\right)_{\Omega}&= (q_{3},w)_{\Omega},\\
\left(K^{-1}\tilde{\bs{F}}_{o},\bs{v}\right)_{\Omega} - \left(P_{o},\nabla \cdot \bs{v}\right)_{\Omega} &= 0.
\end{aligned}
\label{eq:expanded_mixed}
\end{equation}
for all $w\in W_{H,h}$ and $\bs{v} \in V_{ms}$ where $F_o,F_w,F_g$ satisfying the equations in (\ref{eq:expanded_flux}) and  (\ref{eq:numerical_flux}).
\section{Adaptive Numerical Homogenization} \label{sec:adaptive}
In this section, we present an adaptive numerical homogenization scheme for model order reduction or upscaling. We rely upon local numerical homogenization, based on a formal two scale homogenization approach, to evaluate effective properties (permeability, porosity, relative permeability, and capillary pressure) at a coarse scale. We define transient regions as subdomains in space with large gradient in phase saturations to identify the location of the saturation front. A dynamic mesh refinement is used to accurately capture the features of this saturation front directly at the fine scale using an enhanced velocity mixed finite element based domain decomposition approach. This approach allows us to be computationally efficient since pre-processed effective properties, at the coarse scale, are used away from the front while preserving accuracy using fine scale properties along the front.

\subsection{Local Numerical Homogenization}\label{subsec:numhomo}
\begin{figure}[H]
\begin{center}
\includegraphics[width=5cm,trim=0cm 0cm 0cm 0cm, clip]{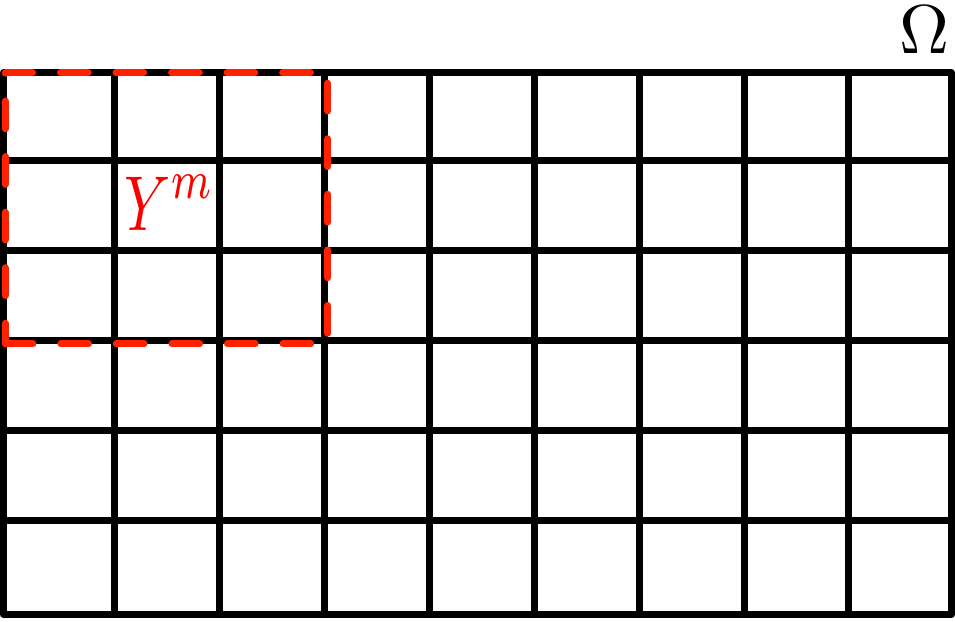} \quad
\includegraphics[width=5cm,trim=0cm 0cm 0cm 0cm, clip]{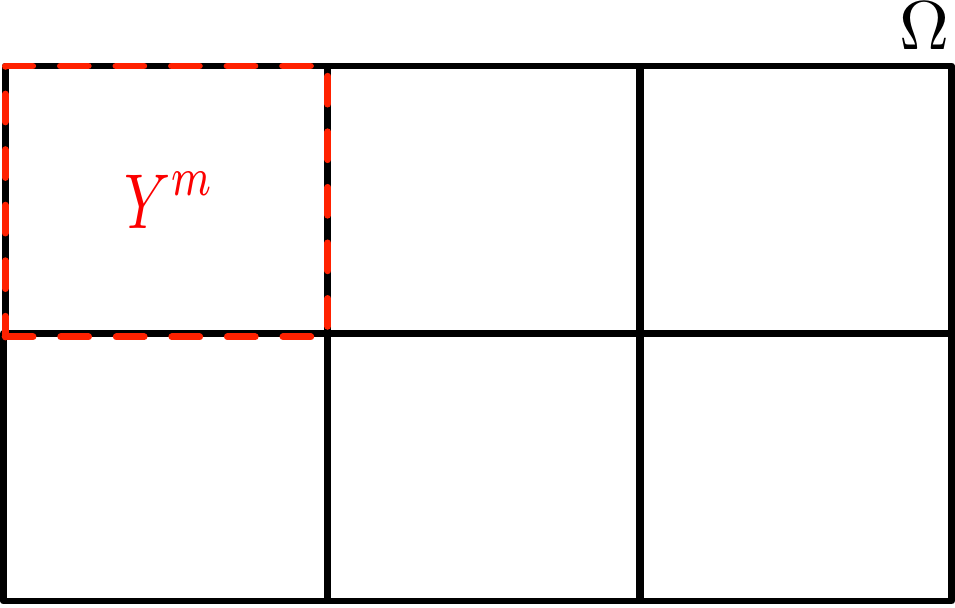}
\caption{Schematic of local numerical homogenization to obtain coarse scale (right) parameters from fine scale (left).}
\label{fig:num_homo}
\end{center}
\end{figure}

Here, we describe the effective equations at the coarse scale for the black oil formulation following the original derivation in \cite{amazianemms,amazianem3as} for two-phase, water-gas type compressible system.

\begin{equation}
\begin{aligned}
\langle\phi\rangle \frac{\partial N^{\text{eff}}_{o}}{\partial t} - \nabla_{x} \cdot \left( \rho_{o}(P^{\text{eff}}_{o})(\Lambda^{\text{eff}}_{o}(\vec{S}^{\text{eff}}_{o})\nabla P^{\text{eff}}_{o1}) \right) = 0\\
\langle\phi\rangle \frac{\partial N^{\text{eff}}_{w}}{\partial t} - \nabla_{x} \cdot \left( \rho_{w}(P^{\text{eff}}_{w})(\Lambda^{\text{eff}}_{w}(\vec{S}^{\text{eff}}_{w})\nabla P^{\text{eff}}_{w}) \right) = 0\\
\langle\phi\rangle \frac{\partial N^{\text{eff}}_{g}}{\partial t} - \nabla_{x} \cdot \left( \rho_{o}(P^{\text{eff}}_{o})\Lambda^{\text{eff}}_{o}(\vec{S}^{\text{eff}}_{o})\nabla P^{\text{eff}}_{o}+ \rho_{g}(P^{\text{eff}}_{g})R_{g}(P^{\text{eff}}_{g})(\Lambda^{\text{eff}}_{g}(\vec{S}^{\text{eff}}_{g})\nabla P^{G}_{g}) \right) = 0\\
\end{aligned}
\end{equation}
Here,
\begin{equation}
\vec{S}_{\alpha} =
\begin{bmatrix}
S_{\alpha,1}\\
\vdots\\
S_{\alpha,k}
\end{bmatrix},
\end{equation}
and $k$ is the number of different rock types with different relative permeability and capillary pressure function descriptions. The subscript $x$ in $\nabla_{x}$ is used to denote the coarse scale operator with the superscript in $A^{\text{eff}}$ to represent the effective properties and coarse scale unknowns. We define global pressures $P^{G}_{k}$ for a rock type $k$ as,
\begin{equation}
\begin{aligned}
P^{\text{eff}}_{\alpha} = P^{G}_{1} + G_{\alpha,1}= \hdots = P^{G}_{k} + G_{\alpha,k}.
\end{aligned}
\end{equation}
Here, $P^{\text{eff}}_{\alpha}$ is the effective pressure at the coarse scale and $P^{G}_{k}$ the global pressure associated with rock type $k$. Further $G_{\alpha,k}$ are defined as,
\begin{equation}
\begin{aligned}
G_{o} (S_{o}) &= G_{o}(0) + \int_{0}^{S_{o}} \frac{\lambda_{g}}{\lambda_{t}} P'_{cgo} d\psi - \int_{0}^{S_{o}} \frac{\lambda_{w}}{\lambda_{t}}P'_{cow} d\psi,\\
G_{w}(S_{o}) &= G_{o}(S_{o}) - P_{cow},\\
 G_{g}(S_{o}) &= G_{o}(S_{o}) + P_{cgo}.
\end{aligned}
\end{equation}
We consider the following auxiliary problems on a subdomain $Y^{m}$ of the larger domain $\Omega = \sum_{m} \cup Y^{m}$ with periodic boundary conditions in the unknowns $\chi_{\alpha}$ as,
\begin{equation}
\begin{aligned}
\nabla_{y} \cdot \left\{ K(y)\lambda_{\alpha}(y,\vec{S}_{\alpha}(y))\left(\mathbb{I} + \nabla_{y} \chi_{\alpha}(y)\right)\right\} &= 0, \quad \chi_{\alpha} \in H^{1}_{per}(Y^{m})\\
\sum_{\alpha} S_{\alpha,k}(y) &= 1
\end{aligned}
\label{eqn:aux}
\end{equation}
Here, $y$ is used to denote fine scale variations in the properties and unknowns. This system of equations can be solved for any $S_{\alpha,k} \in [S^{\text{irr}}_{\alpha,k},S^{\text{res}}_{\alpha,k}]$ between irreducible and residual phase saturations for the unknowns $\chi_{\alpha}$. The effective tensor $\Lambda^{\text{eff}}_{\alpha}$ is then evaluated using,
\begin{equation}
\Lambda^{\text{eff}}_{\alpha}(\vec{S}_{\alpha}(y)) = \int_{Y^{m}} K(y) \lambda_{\alpha}(y,\vec{S}_{\alpha}(y))\left\{\mathbb{I}+\nabla_{y} \chi_{\alpha}(y)\right\}dy.
\end{equation}
The overall effective porosity $\langle \phi \rangle$ is given by,
\begin{equation}
\begin{aligned}
\langle \phi \rangle = \frac{1}{|Y^{m}|}\int_{Y^{m}}\phi(y)dy.
\end{aligned}
\end{equation}
We define effective porosity for the rock type $k$ ($\bar{\phi}_{k}$) in subdomain $Y^{m}$ as,
\begin{equation}
\begin{aligned}
\bar{\phi}_{k} = \frac{1}{|Y^{m}_{k}|}\int_{Y^{m}_{k}}\phi(y)dy.
\end{aligned}
\end{equation}
Here, $Y^{m}_{k}$ is the subdomain of $Y^{m}$ that belongs to rock type denoted by index $k$. Further, $S^{\text{eff}}_{\alpha}$ is defined as,
\begin{equation}
S_{\alpha}^{\text{eff}} = \sum_{k} \frac{|Y^{m}_{k}|}{\langle \phi \rangle}\bar{\phi}_{k} S_{\alpha,k}.
\label{eqn:seff}
\end{equation}
The effective concentration $N_{\alpha}^{\text{eff}}$ are related to the effective saturations $S_{\alpha}^{\text{eff}}$ using Eqn. \eqref{eqn:consatrel}. The phase saturations in different rock types $S_{\alpha,k}$ are related using,
\begin{equation}
\begin{aligned}
P_{cow,1}(S_{o,1}) = \hdots = P_{cow,k}(S_{o,1}),\\
P_{cgo,1}(S_{o,1}) = \hdots = P_{cgo,k}(S_{o,1}).
\end{aligned}
\label{eqn:cappresk}
\end{equation}
Equations \eqref{eqn:seff} and \eqref{eqn:cappresk} are used to eliminate multiple phase saturations, associated with different rock types, in favor of one effective phase saturation. The effective capillary pressures can be evaluated as,
\begin{equation}
\begin{aligned}
P^{\text{eff}}_{cow}(S^{\text{eff}}_{o}) = G_{o,1}(S^{\text{eff}}_{o}) - G_{w,1}(S^{\text{eff}}_{o}) = \hdots = G_{o,k}(S^{\text{eff}}_{o}) - G_{w,k}(S^{\text{eff}}_{o}),\\
P^{\text{eff}}_{cgo}(S^{\text{eff}}_{o}) = G_{g,1}(S^{\text{eff}}_{o}) - G_{o,1}(S^{\text{eff}}_{o}) = \hdots = G_{g,k}(S^{\text{eff}}_{o}) - G_{o,k}(S^{\text{eff}}_{o}).
\end{aligned}
\end{equation}
We solve the aforementioned auxiliary problems Eqns. \eqref{eqn:aux} for the non-overlapping subdomains $Y^{m}$ as shown in Figure \ref{fig:num_homo}.
The effective properties at the coarse scale namely: capillary pressure $P^{\text{eff}}_{c}$, porosity $\langle \phi \rangle$, and mobility tensor $\Lambda^{\text{eff}}_{\alpha}$, are computed only once as a pre-processor step prior to the numerical simulation.
\subsection{Adaptivity Criterion} \label{subsec:criterion}
We describe an indicator function to track the location of the saturation front in order to later perform a domain decomposition into transient and non-transient regions.  An ad hoc criterion for identifying the location of the saturation front can be easily defined using a gradient in saturation between saturation at a given point in space and its nearest neighbor. One such criterion is using a maximum of absolute of difference between a saturation $S$ at an element and its adjacent elements at the previous time step $n$. We define the neighboring elements collectively as, $\Omega_{neighbor}(x)=\{y: y\in E_j, |\partial E_i \cap \partial E_j |   \neq \emptyset, \text{if } x\in E_i  \}$. Then the adaptivity criteria can be written as,
			\begin{align} \label{eq:criteria2}
			\textcolor{blue}{\Omega_f}=\left\{\textbf{x}: \text{max} |S_{w}^{n}(x)-S_{w}^{n}(y)|>\epsilon_{adap} \quad \forall y \in \Omega_{neighbor}(x)  \right\}
			\end{align}
Here, $E_{i}$ and $E_{j}$ represent an element and its neighbors with $\epsilon_{adap}$ as the threshold value above which a domain is marked as a transient region. Please note that this type of adaptivity criterion has been used by others \cite{aarnes2006adaptive} to reduce computational costs in a similar sense.

\subsection{Adaptive Mesh Refinement}
\begin{figure}[H]
\begin{center}
\includegraphics[width=6cm,trim=0cm 0cm 0cm 0cm, clip]{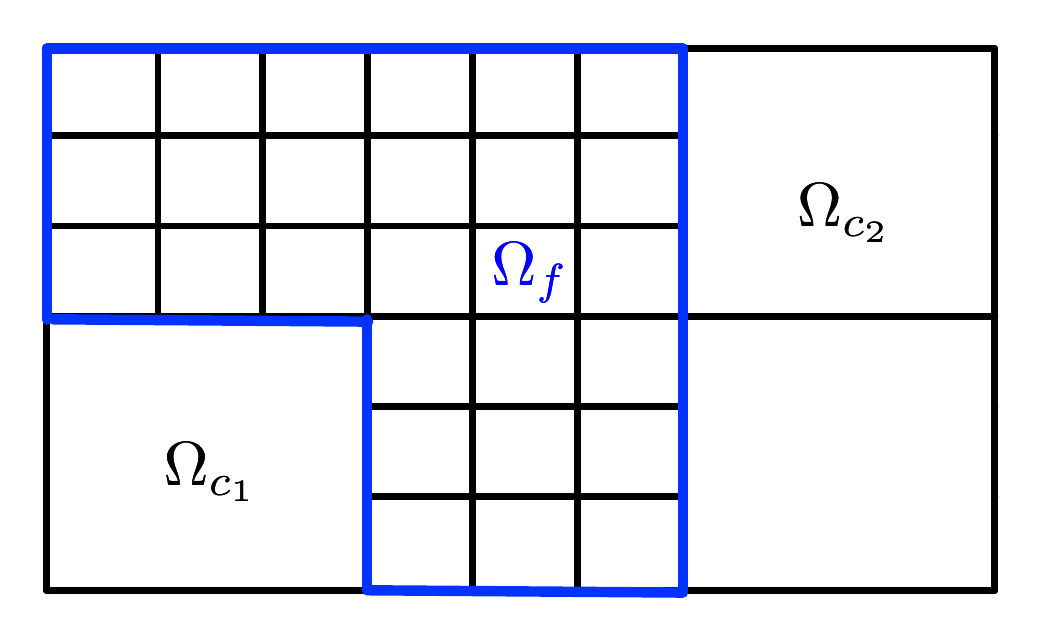}
\caption{Schematic of adaptive mesh refinement with coarse ($\Omega_{c}$) and fine ($\textcolor{blue}{\Omega_{f}}$) domains.}
\label{fig:adap}
\end{center}
\end{figure}

Based upon the above criteria we divide the domain ($\Omega$) into non-overlapping, transient ($\textcolor{blue}{\Omega_{f}}$) and non-transient ($\Omega_{c}$) subdomains to solve flow and transport problems at the fine and coarse scales, respectively. Figure \ref{fig:adap} shows a schematic of domain decomposition approach into fine (\textcolor{blue}{$\Omega_{f}$}) and coarse ($\Omega_{c}$) subdomains. In what follows, coarse and non-transient, and fine and transient can be used interchangeably to refer to a subdomain. The coarse and fine subdomain problems are then coupled at the interface using the enhanced velocity mixed finite element (EVMFE) method  spatial discretization described in \cite{wheeler2002enhanced}. This multiblock, domain decomposition approach is strongly mass conservative at the interface between fine and coarse domains and hence preserves local mass conservation property of the mixed finite element scheme. The EVMFE scheme has been used previously for a number of fluid flow and transport problems \cite{thomas2011enhanced} including equation of state (EOS) based compositional flow. As mentioned before, in the fine domain we fine scale properties directly while relying upon effective properties in the coarse domain obtained from local numerical homogenization described in subsection \ref{subsec:numhomo}.

\section{Generalized Mixed Multiscale Method}\label{sec:GMMM}
In this section, we will introduce the generalized mixed multiscale method for the black oil model. The main idea of the method is to construct a coarse-grid multiscale space which can approximate the solution of the black oil model accurately.
Before presenting the method, we need to introduce the some notations of the coarse grid and the fine grid first. We let $\mathcal{T}_{H}$ be a partition
of the computational domain $\Omega$ into finite elements. $H$ is the mesh size of $\mathcal{T}_H$. We call this partition as
coarse grid. The fine grid partition is then defined as a refinement of the coarse grid. We denote the fine grid partition by $\mathcal{T}_{h}$. Similarly, $h$ is the fine mesh size of $\mathcal{T}_h $. Next, we will denote $\mathcal{E}_H $ as the set of all coarse edges. The set of fine edges is then denoted as $\mathcal{E}_h$.

In problem (\ref{eq:expanded_mixed}), there are four finite element spaces which are used to approximate the oil pressure, water saturation, gas saturation, and the pseudo-velocity with respectively.
In this method, we construct the multiscale space for approximating the pseudo-velocity. For the pressure space and the saturation space, we will use the piecewise constant space corresponding to the coarse-grid. Since the saturations satisfy the transport equations with capillary pressure effects, it is not easy to construct an multiscale space for the saturations on the coarse grid. To approximate the saturations accurately, we will use an adaptive approach to define the finite element space for pressure and saturations. We will use a fine grid piecewise constant space in some of the regions and use coarse grid piecewise constant space in the rest of the region. We will discuss the detail of the adaptive method in subsection \ref{sec:adaptive_update}. Once we constructed the space for the saturations, pressure and velocity, the resulting method is then described as; find $(P_{o},S_{w},S_{g},\bs{\tilde{F}}_{o}) \in W_{H,h}\times W_{H,h}\times W_{H,h}\times V_{H,h}$ such that

\[
(P_{o},S_{w},S_{g},\bs{\tilde{F}})  \text{ is a solution of the equations in (\ref{eq:expanded_mixed})}
\]
for all $(w,v)\in W_{H,h}\times V_{H,h}$.
$W_{H,h}$ denotes the space of all piecewise constant functions corresponding to the adaptive mesh, that is
\[
W_{H,h}=\{q\in L^{2}(\Omega)|\;q|_{K}\in P_{0}(K),\;\text{for} K\in\mathcal{T}_{h}(\Omega_h)\cup \mathcal{T}_{H}(\Omega\backslash\Omega_h) \}.
\]
where $\mathcal{T}_h(\Omega_h) = \{K\in \mathcal{T}_h|\;K\subset \Omega_h\}$ and $\mathcal{T}_{H}(\Omega\backslash\Omega_h)= \{K\in \mathcal{T_H}|\;K\subset \Omega\backslash\Omega_h\}$ and 
$V_{H,h}$ denotes the sum of the multiscale velocity space and fine grid velocity space corresponding to the fine grid on $\Omega_h$, that is,
\[
V_{H,h}=V_{ms}+\{v\in H(\text{div},\Omega)|\;v|_{f}\in RT_{0}(f),\text{ for } f\in\mathcal{T}_{h}(\Omega_h)\text{ and } v\cdot n|_{\partial F} = 0 \; \forall F\in \mathcal{T}_H\}.
\]
We remark that it is possible to increase the accuracy of the approximation
for the pressure solution if we enrich the space by multsicale basis functions in coarse grid level
which follow the idea in \cite{chung2018multiscale}. To keep the proposed method simple,
we will not consider the enriched space in this paper.

The most important procedure of this method is the construction of
the multiscale velocity space. The multiscale velocity space is constructed
in two steps. The first step is constructing a snapshot space for
the velocity field. The second step is preforming a dimension reduction
process on the snapshot space to construct the multiscale velocity
space. The snapshot space contain functions which satisfy the static
flow equation locally with different flux boundary conditions on coarse
edges. The multiscale velocity space is constructed by performing
local spectral decomposition on the snapshot space to select the dominant
eigenfunctions. We will discuss the detail of the construction process
in Section \ref{sec:basis_construction}.

After obtaining the mutiscale space, we will solve the problem (\ref{eq:expanded_mixed})
by using this space.
To improve the computational speed of the method, one can solve the pressure, the saturations and the velocity in the coarse grid space first. That is, we consider the fine grid regions to be empty ($\Omega_h=\emptyset$). Since the size of the coarse grid problem is small, we can get a coarse grid approximation of the solution quickly. Then we can use the coarse grid solution to be the initial guess for the adaptive grid solver. Moreover, the problem in adaptive grid can be solved by using the coarse grid system as a preconditioner which give us a efficient solver for the adaptive grid problem.

\subsection{Multiscale basis construction}
\label{sec:basis_construction}
In this section, we will discuss the construction of the multiscale
basis functions for the velocity space. We will first introduce the
snapshot space, which contains a set of basis function satisfying the
local static flow equation with some boundary conditions on the coarse
grid edge. Next, we will introduce a dimension reduction technique
to construct a low dimension local multiscale space which is the span
of some dominant modes in the local snapshot space. These dominant
modes are selected by solving a local spectral problem for each local
domain.

\subsubsection*{Snapshot space}

In this subsection, we will discuss the construction of the snapshot
space. In this paper, we will introduce two choices of the snapshot
space. One of the snapshot space contains the solutions of local static
flow problems with all possible boundary conditions on the coarse
edge up to the fine-grid resolution. The another snapshot space is
a subspace of the the first choice of snapshot space. Instead of containing
the basis function with all of possible boundary conditions, it only
contains the basis functions with certain boundary condition on the
coarse edge. To select these boundary condition, we will solve the
local problem in an oversampling region with some random boundary
conditions. The boundary conditions for the snapshot functions are
the restriction of the oversampling solutions on the coarse edge.
Next, we will explain the construction of the first choice of snapshot
space in detail.

Let $E_{i}\in\mathcal{E}_{H}$ be a coarse edge. The coarse edge $E_{i}$
can be written as a union of some fine edge, that is, $E_{i}=\cup_{j=1}^{J_{i}}e_{j}$
where $e_{j}$ are fine edge in $\mathcal{E}_{h}$ and $J_{i}$ is
the number of fine edge on $E_{i}$. We then denote the coarse neighborhood
$\omega_{i}$ as $\omega_{i}=\cup\{F\in\mathcal{T}_{H}|\bar{F}\cap E_{i}\neq\emptyset\}$.
We find $(\psi_{j}^{(i)},q_{j}^{(i)})$ by solving the following problem
on the coarse neighborhood $\omega_{i}$
\begin{equation}\label{eq:local_problem}
\begin{array}{rll}
K^{-1}\psi_{j}^{(i)}+\nabla q_{j}^{(i)} & =0 & \text{ in }F\in\mathcal{T}_{H},\text{ with }F\subset\omega_{i}\\
\nabla\cdot(\psi_{j}^{(i)}) & =\alpha_{j}^{(i)} & \text{ in }F\in\mathcal{T}_{H},\text{ with }F\subset\omega_{i}
\end{array}
\end{equation}
subject to the boundary condition

\[
\begin{array}{rll}
\psi_{j}^{(i)}\cdot n_{\partial\omega_{i}} & =0 & \text{ on }\partial\omega_{i}\\
\psi_{j}^{(i)}\cdot n_{E_{i}} & =\delta_{j}^{(i)} & \text{ in }E_{i}
\end{array}
\]
where $n_{\partial\omega_{i}}$ is the unit outward normal on $\partial\omega_{i}$
and $n_{E_{i}}$ is a fixed unit normal on $E_{i}$. The function
$\delta_{j}^{(i)}$ is a piecewise constant function with respect
to the fine-grid on $E_{i}$ with value $1$ on $e_{j}$ and value
$0$ on other fine edge, namely,
\[
\delta_{j}^{(i)}=\begin{cases}
1, & \text{ on }e_{j},\\
0, & \text{ on }e_{l}\neq e_{j}
\end{cases}.
\]
The function $\alpha_{j}^{(i)}$ is piecewise constant function with
respect to the coarse-grid with
\[
\alpha_{j}^{(i)}|_{F}=n_{E_{i}}\cdot n_{\partial F}\cfrac{|e_{j}|}{|F|}\;\text{ in }F\in\mathcal{T}_{H},\text{ with }F\subset\omega_{i}.
\]
 Since the snapshot basis function $\psi_{j}^{(i)}$ satisfies $\psi_{j}^{(i)}\cdot n_{\partial\omega_{i}}=0$
on $\partial\omega_{i}$ , $\psi_{j}^{(i)}$ can extend to $\Omega\backslash\omega_{i}$
by zero and the extension is still in $H(\text{div},\Omega)$. The
local snapshot space $V_{snap}^{(i)}$ is then defined as
\[
V_{snap}^{(i)}=\text{span}\{\psi_{j}^{(i)}|\;1\leq j\leq J_{i}\}.
\]
The snapshot space $V_{snap}$ is defined as
\[
V_{snap}=\oplus V_{snap}^{(i)}=\text{span}\{\psi_{j}^{(i)}|\;1\leq j\leq J_{i},\;1\leq i\leq N_{e}\}.
\]

\begin{figure}[htb]
\centering
\includegraphics[scale=0.5]{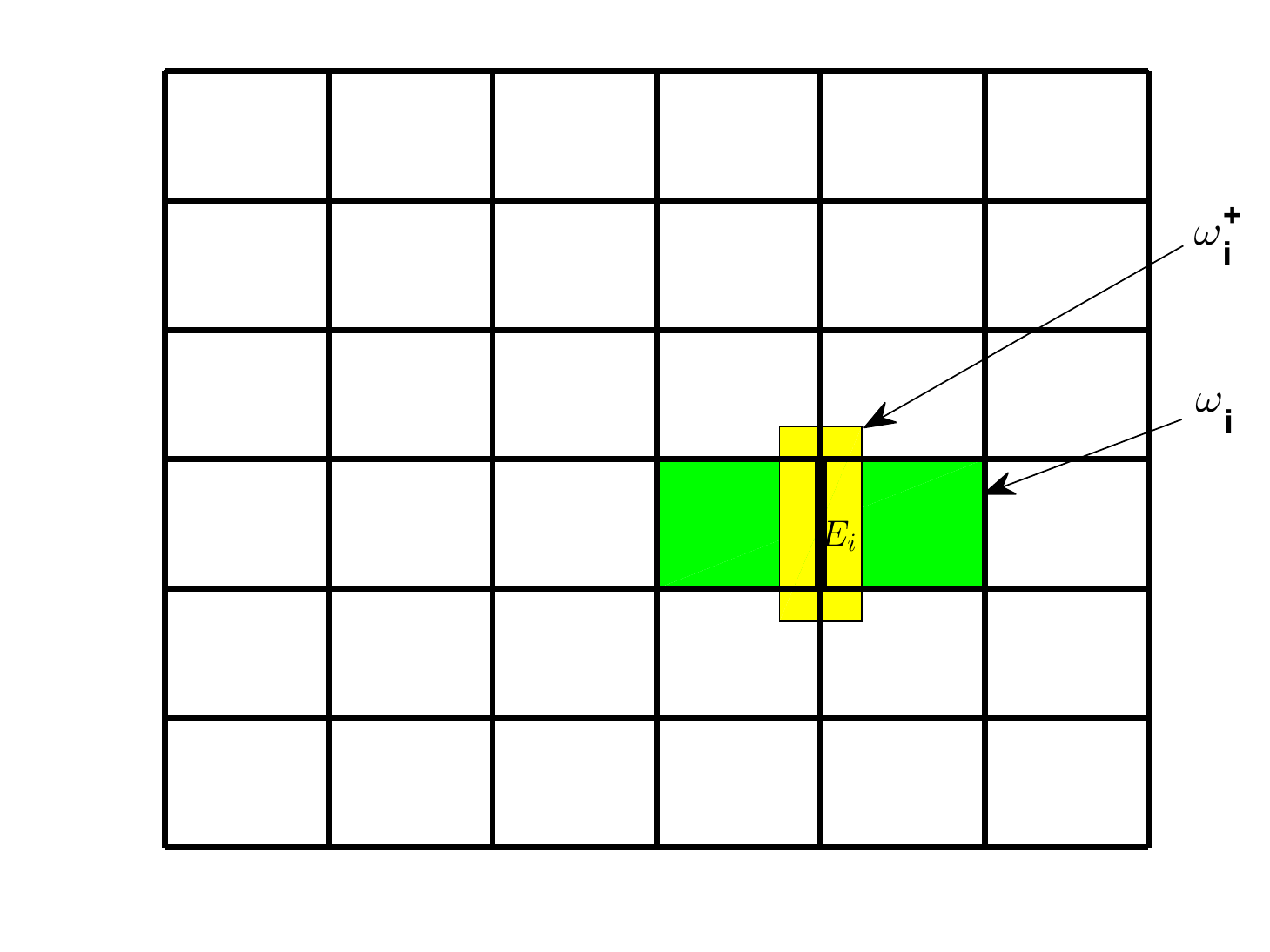}
\label{fig:illustration_over}

\caption{An illustration of the oversampling domain}
\end{figure}

Next, we will discuss the second choice of the snapshot space. For
a coarse edge $E_{i}$, we can define a oversampling domain $\omega_{i}^{+}\supset E_{i}$
with $\omega_{i}^{+}$ is an union of the fine-grid element, that
is, $\omega_{i}^{+}=\cup_{F\in I}K$ for some $I\subset\mathcal{T}_{h}$.
We will show an illustration of the oversampling domain in figure
\ref{fig:illustration_over}.

We find $(\psi_{j}^{(i),+},q_{j}^{(i),+})$ by solving the following
problem on the oversampling domain $\omega_{i}^{+}$
\[
\begin{array}{rll}
K^{-1}\psi_{j}^{(i),+}+\nabla q_{j}^{(i),+} & =0 & \text{ in }\omega_{i}^{+}\\
\nabla\cdot(\psi_{j}^{(i),+}) & =\beta_{j}^{(i)} & \text{ in }\omega_{i}^{+}
\end{array}
\]
with boundary condition $\psi_{j}^{(i),+}\cdot n_{\partial\omega_{i}^{+}}=r_{j}^{(i)}$
for $1\leq j\leq J_{i}^{+}$ where $r_{j}^{(i)}$ is a piecewise constant
with respect to the fine-grid on $\partial\omega_{i}^{+}$ with Gaussian
random value and $\beta_{j}^{(i)}=\cfrac{\int_{\partial\omega_{i}^{+}}r_{j}^{(i)}}{|\omega_{i}^{+}|}$.
The snapshot basis function $\psi_{j}^{(i)}$ is then the solution
of the local problem (\ref{eq:local_problem}) with boundary condition

\[
\begin{array}{rll}
\psi_{j}^{(i)}\cdot n_{\partial\omega_{i}} & =0 & \text{ on }\partial\omega_{i}\\
\psi_{j}^{(i)}\cdot n_{E_{i}} & =\psi_{j}^{(i),+}\cdot n_{E_{i}} & \text{ in }E_{i}.
\end{array}
\]
Similarly, the local snapshot space $V_{snap}^{(i)}$ is defined as
the span of the snapshot basis function $\psi_{j}^{(i)}$ and the
snapshot space $V_{snap}$ is defined as the sum of the local snapshot
space $V_{snap}^{(i)}$.

\subsubsection*{Multiscale space}

In this section, we will discuss the construction of the multiscale
space $V_{ms}$. Following the framework of \cite{chung2015mixed}, we will perform
a dimension reduction on the snapshot space by using some local spectral
problems. The spectral problems are used to select some dominant modes
from the local snapshot spaces $V_{snap}^{(i)}$. These local dominant
modes span the local multiscale space and the multiscale is the sum
of these local multiscale space. For each coarse edge $E_{i}$, we
will solve a local spectral problem which find $(\phi_{j}^{(i)},\lambda_{j}^{(i)})\in V_{snap}^{(i)}\times\mathbb{R}$
such that

\[
a^{(i)}(\phi_{j}^{(i)},v)=\lambda_{j}^{(i)}s^{(i)}(\phi_{j}^{(i)},v)\;\forall v\in V_{snap}^{(i)},
\]
 where $a^{(i)}(\cdot,\cdot)$ and $s^{(i)}(\cdot,\cdot)$ are two
symmetric positive definite bilinear forms. There are several choices
for the local spectral problem. In \cite{chung2015mixed, chan2016adaptive}, authors suggested
different spectral problem. In this paper, we will select one of these
spectral problems. We consider the bilinear forms $a^{(i)}:V_{snap}^{(i)}\times V_{snap}^{(i)}\rightarrow\mathbb{R}$
and $s^{(i)}:V_{snap}^{(i)}\times V_{snap}^{(i)}\rightarrow\mathbb{R}$
are defined as

\[
\begin{array}{rll}
a^{(i)}(u,v)=\int_{E_{i}}k^{-1}(u\cdot n_{E_{i}})(v\cdot n_{E_{i}}) & \text{ and } & s^{(i)}(u,v)=\int_{\omega_{i}}k^{-1}u\cdot v+\int_{\omega_{i}}(\nabla\cdot u)(\nabla\cdot v)\end{array}.
\]
We assume the eigenvalue are arranged in ascending order, that is,
\[
0<\lambda_{1}^{(i)}\leq\lambda_{2}^{(i)}\leq\lambda_{3}^{(i)}\leq\cdots.
\]
We will use the first $L_{i}$ eigenfunctions to span our local multiscale
space $V_{ms}^{(i)}$, namely,
\[
V_{ms}^{(i)}=\text{span}\{\phi_{j}^{(i)}|\;1\leq j\leq L_{j}\},
\]
where $L_{i}$ is the dimension of the local multsicale space $V_{ms}^{(i)}$.
The multiscale space $V_{ms}$ is then defined as the sum of the local
multiscale space $V_{ms}^{(i)}$, that is,
\[
V_{ms}=\oplus_{i}V_{ms}^{(i)},
\]
and the total dimension of $V_{ms}$ is the sum of $L_{i}$, namely,
$\text{dim}(V_{ms})=\sum_{i}L_{i}$. In practice, we normally use
a small number of basis functions. There are two ways to determine
how many basis functions are required. One is using the eigenvalue
to estimate the number of basis functions $L_{i}^{(i)}$. Following
the analysis in \cite{chung2015mixed}, the approximation error of the multsicale
space is related to the minimum of the eigenvalue $\lambda_{L_{i}+1}^{(i)}$.
If $\min_{i}\{\lambda_{L_{i}+1}^{(i)}\}$ is small, the multiscale
space may give an inaccurate approximation of the solution. Therefore,
we will need to use sufficiently many multiscale basis functions to
ensure $\min\{\lambda_{L_{i}+1}^{(i)}\}$ does not closed to zero.
Another way to determine the number $L_{i}$ is using a posteriori
error estimate to locate which region requires more basis functions.
We can apply this approach in both the offline stage and the online
stage. In the offline stage, we can solve a number of sample problems
in a short time period and use the numerical result to preform the
posteriori error estimate. In the online stage, we can use the online
numerical result to modify the number of basis functions we used.

\subsection{Adaptive approach}
\label{sec:adaptive_update}
In this section, we will present an adaptive griding
approach to construct a finite element space which gives a good approximation for the saturations and the pressure. We will determent the fine grid region $\Omega_h$ by using a coarse grid indicator. There are different choices for the indicator. One of the indicators is describe in subsection \ref{subsec:criterion}. In this subsection, we introduce a residual based indicator to locate the fine grid region. We will use a fine grid finite element space in the region with large residual. That is
\[
\Omega_h = \cup\{K\in \mathcal{T}_H|\;\|R_K\|>\theta(\max_K\{\|R_K\|\})\},
\] where the residual operator $R_{\alpha}$ is defined as
\[
R_{K,\alpha}=\left(\frac{\partial}{\partial t} \left(\phi N_{\alpha}(S_{w,ms},S_{g,ms},P_{o,ms}) \right) + \nabla \cdot \bs{F}_{\alpha,ms} - q_{\alpha},w\right)_{\Omega}
\]
and the residual norm is defined as
\[
\|R_{K}\|=\max_{\alpha}\{\sup_{w\in W_{h}(K)}\cfrac{|R_{\alpha}(w)|}{\|w\|_{L^{2}(K)}}\}.
\]
where $(S_{\alpha,ms},P_{o,ms},\bs{F}_{\alpha,ms})$ is the solution of the problem (\ref{eq:expanded_mixed}) in coarse grid space.
We remark that these residual operators can also be used to construct online multiscale basis functions to enrich the coarse grid space.
%
%

\section{Numerical Results}\label{sec:results}

\begin{figure}[H]
\begin{center}
\includegraphics[width=7.5cm,trim=2.5cm 8cm 2.5cm 7cm, clip]{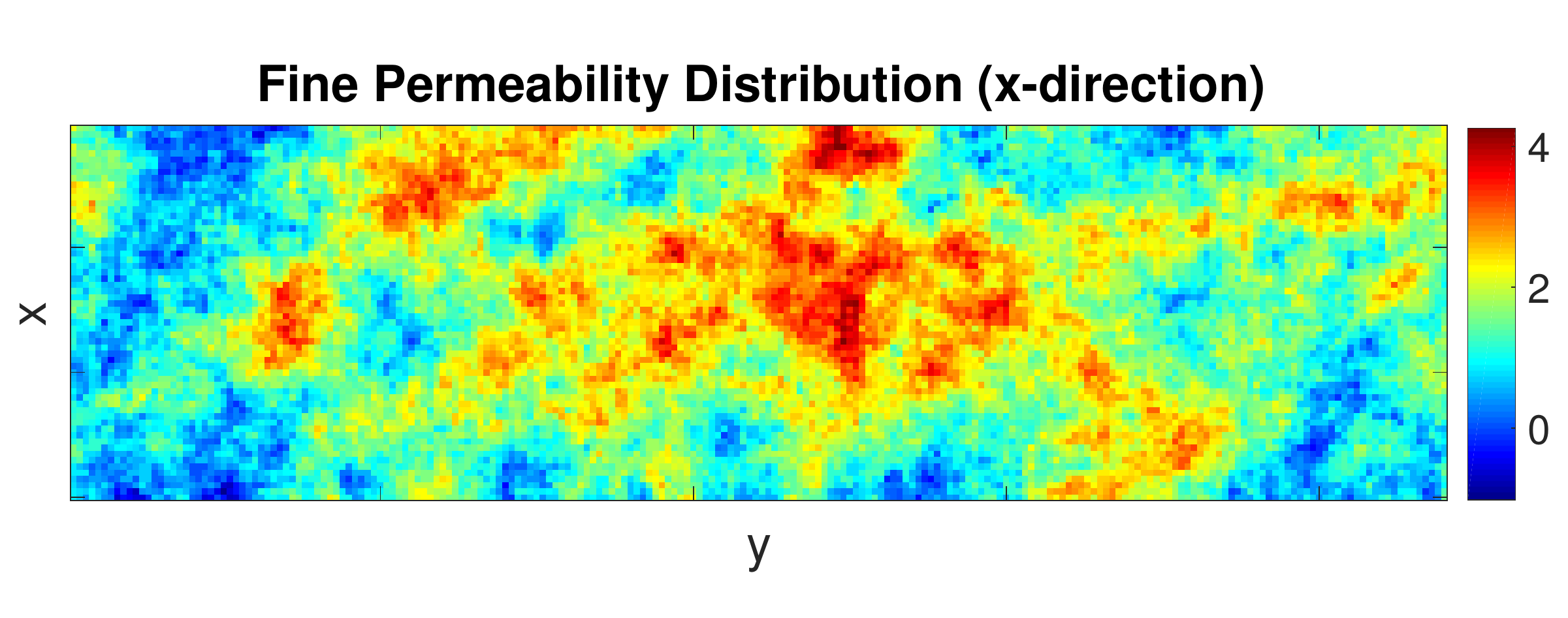}
\includegraphics[width=7.5cm,trim=2.5cm 8cm 2.5cm 7cm, clip]{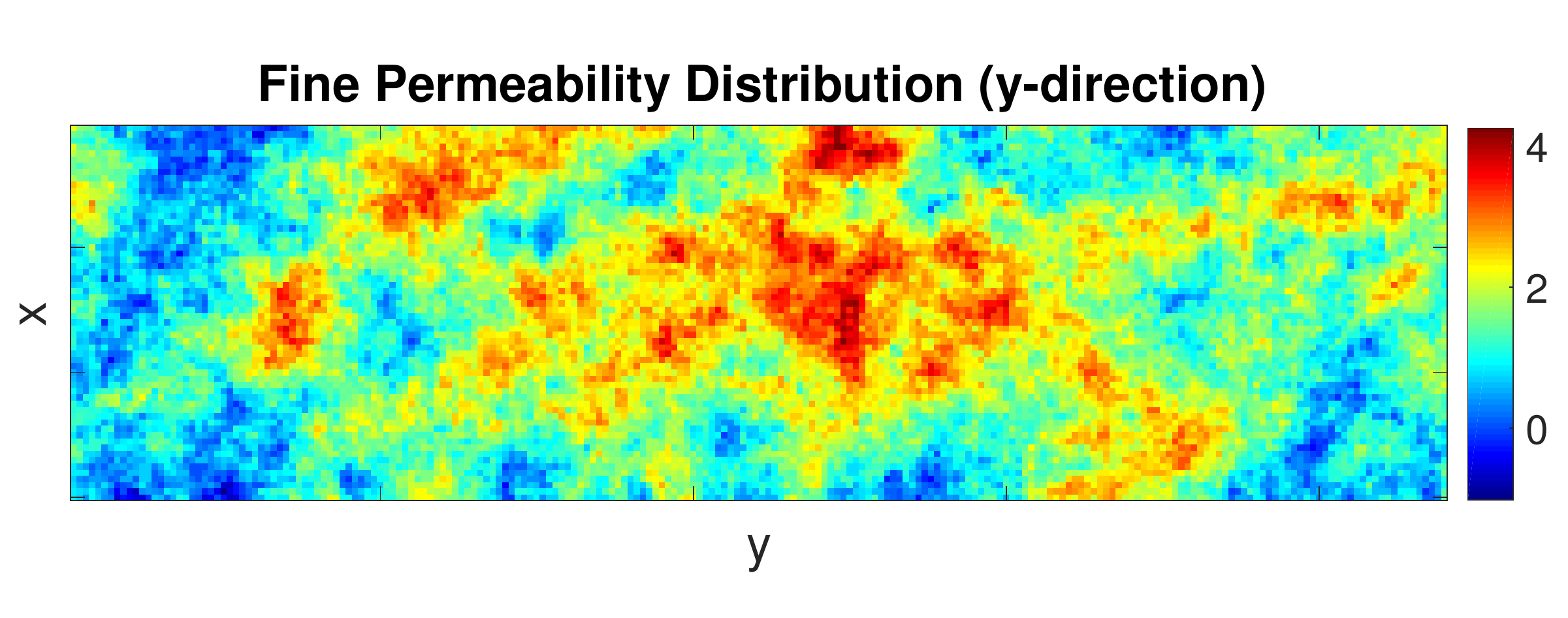}
\caption{Permeability distribution ($log_{10}$ scale) from layer 20 of the SPE10 dataset.}
\label{fig:permdist}
\end{center}
\end{figure}

\begin{figure}[H]
\begin{center}
\includegraphics[width=7.cm,trim=0cm 0cm 0cm 0cm, clip]{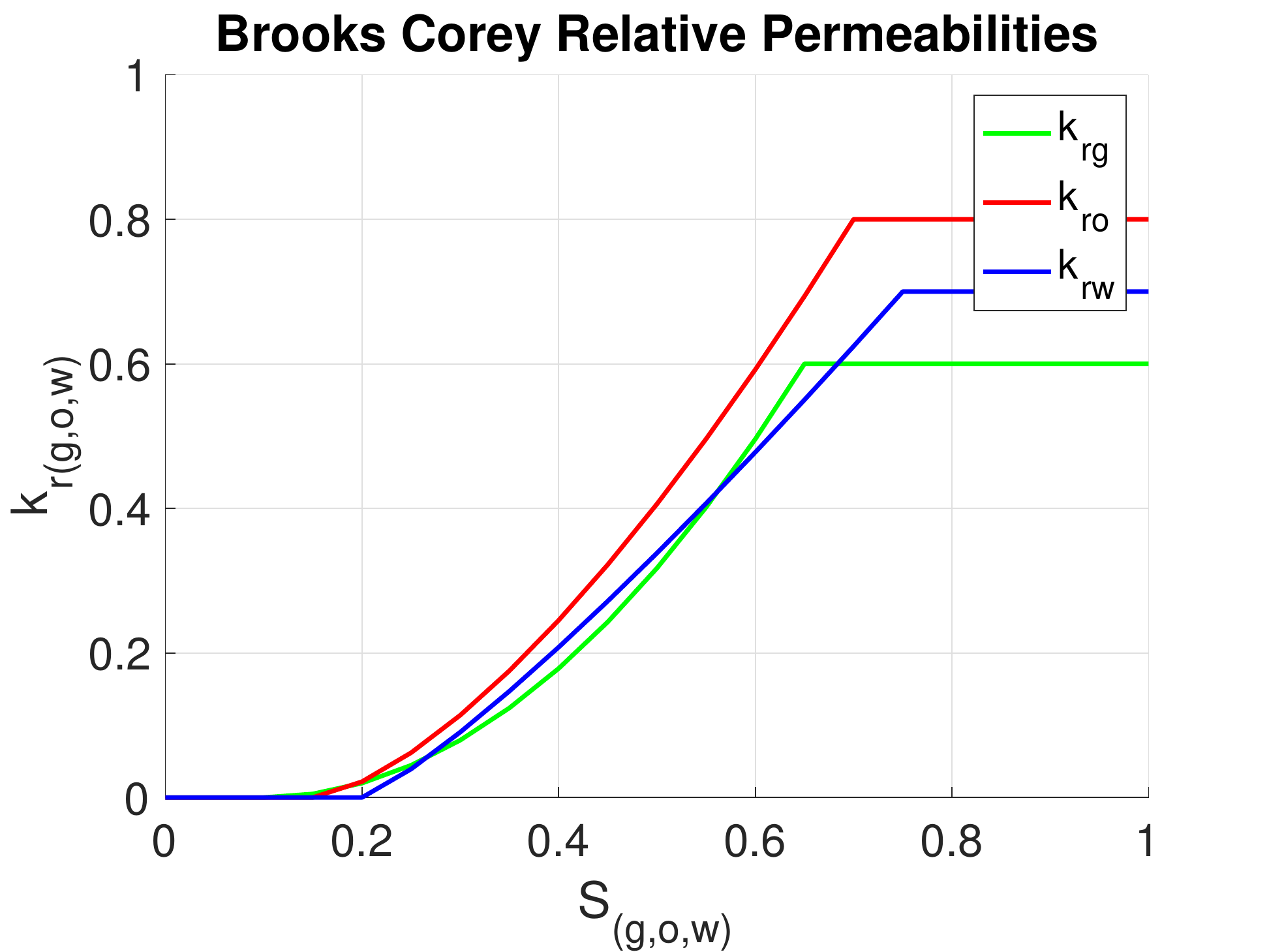}
\includegraphics[width=7.cm,trim=0cm 0cm 0cm 0cm, clip]{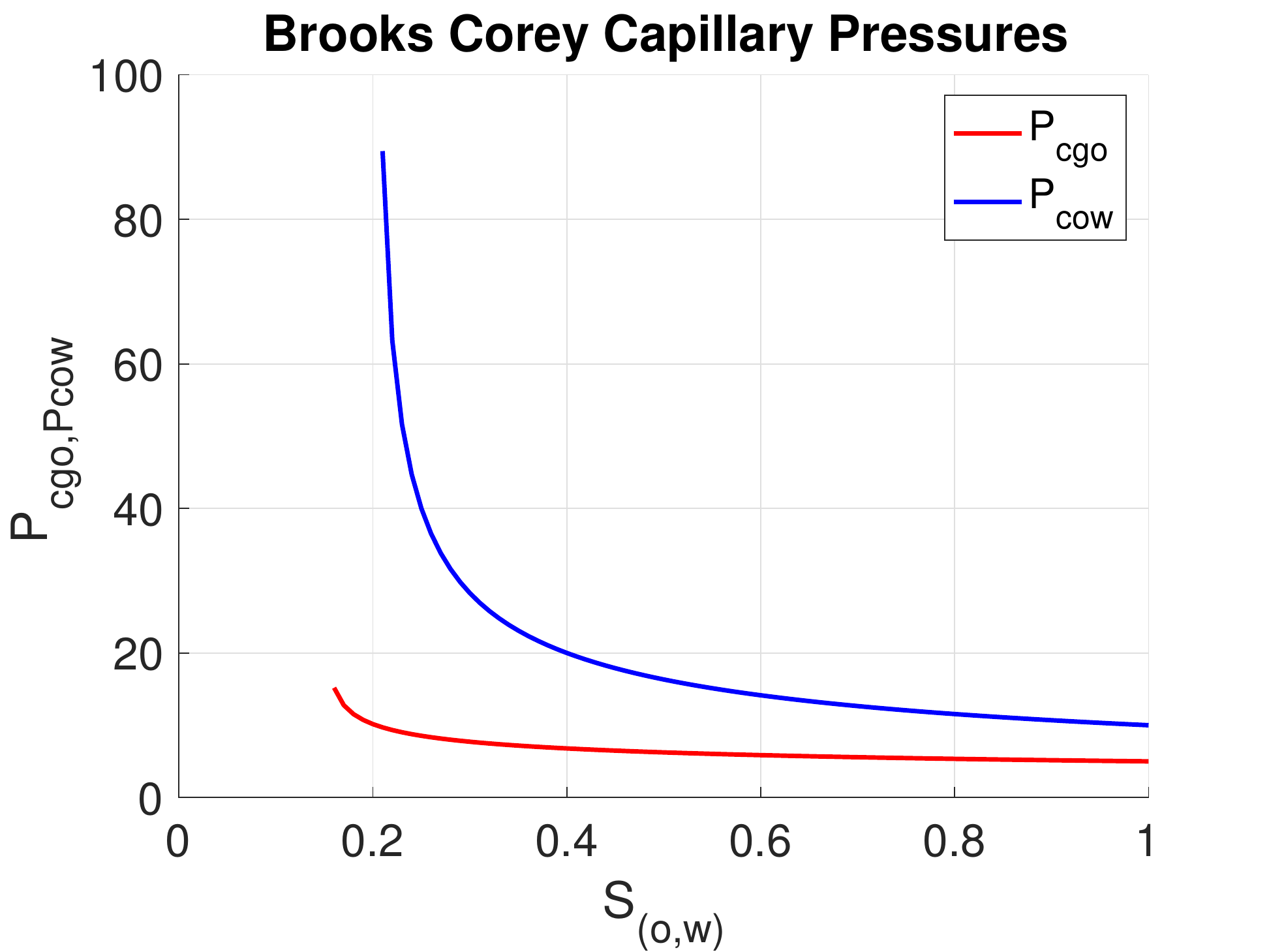}
\caption{Brooks-Corey relative permeability and capillary pressure curves.}
\label{fig:relcap}
\end{center}
\end{figure}

\begin{figure}[H]
\begin{center}
\includegraphics[width=7.cm,trim=1cm 0cm 1cm 0cm, clip]{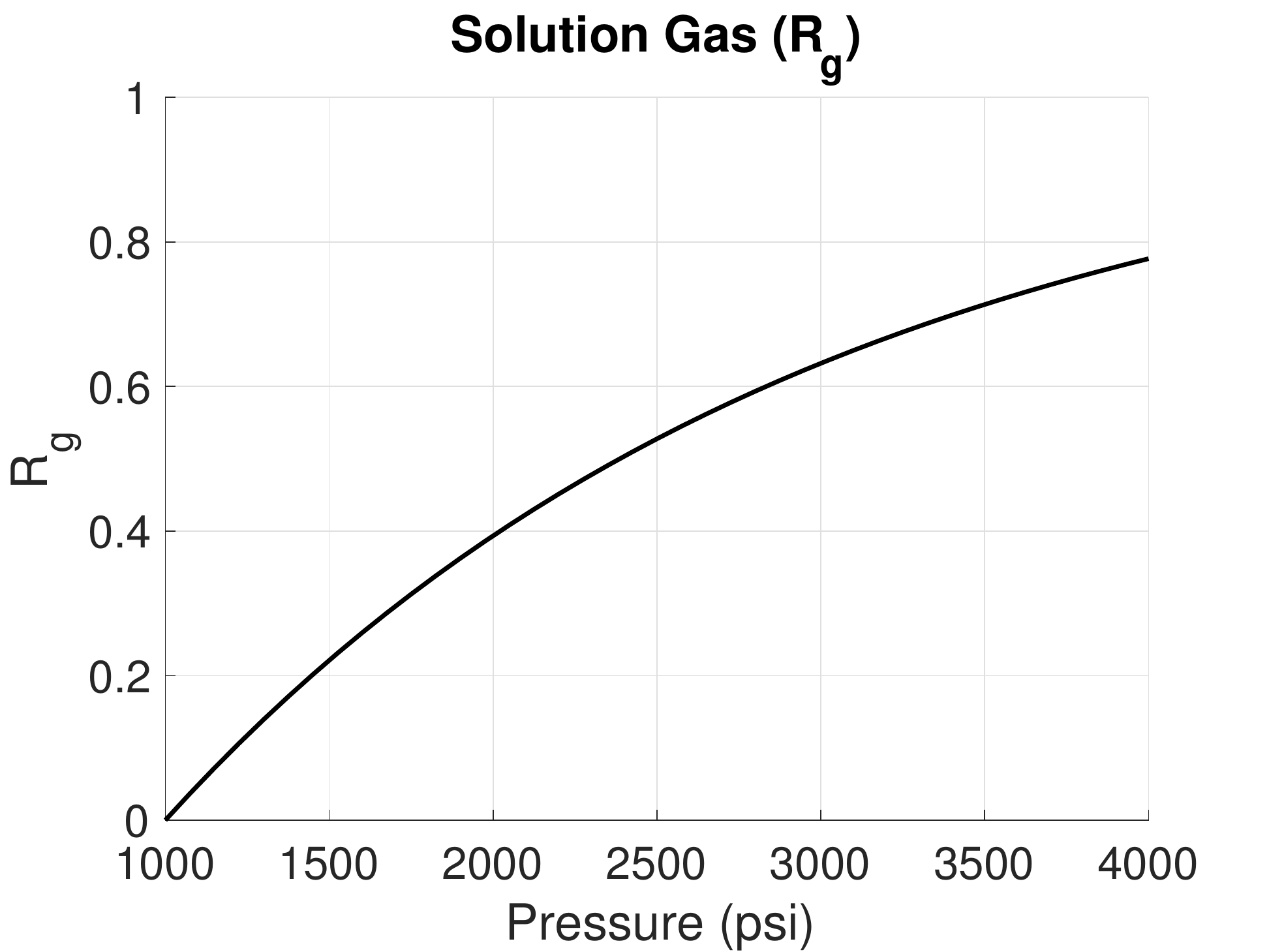}
\caption{Solution gas $R_{g}$ variation with pressure.}
\label{fig:solgas}
\end{center}
\end{figure}

In this section, we present numerical experiments for the adaptive homogenization and generalized multiscale approaches using an augmented dataset from the 10$^\text{th}$ SPE comparative project \cite{spe10}. Figure \ref{fig:permdist} shows the spatial distribution of x and y direction, diagonal components of the absolute permeability tensor. The reservoir domain is 120ft$\times$30ft with coarse and fine scale grid discretizations of 22$\times$6 and 220$\times$60, respectively. The coarse and fine grid elements are consequently 5ft$\times$5ft and 0.5ft$\times$0.5ft. Although not restrictive, the reservoir porosity is assumed to be homogeneous with a value of 0.2.
The gas phase is assumed to be fully compressible with slightly compressible description for the oil and water phases described by equations,
\begin{equation}
\begin{aligned}
\rho_{g} & = c_{g}P_{g},\\
\rho_{o} & = \rho_{o}^{0} exp^{c_{o}P_{o}},\\
\rho_{w} & = \rho_{w}^{0}exp^{c_{w}P_{w}}.
\end{aligned}
\label{eqn:density}
\end{equation}
Here, the gas, oil, and water phase compressibility is taken to be $c_{g}$ = 4.7$\times$ 10$^-3$, $c_{o}$ = 1$\times$10$^{-4}$, and $c_{w}$ = 1$\times$10$^{-6}~$psi$^{-1}$ , respectively. Further, the oil and water phase density at the standard conditions are $\rho_{o}^{0}$ = 56.0 and $\rho_{w}^{0}$ = 66.5 lbs/cuft, respectively. The fluid viscosities are assumed to be constant at $\mu_{g}$ = 0.018, $\mu_{o}$ = 2, and $\mu_{w}$ = 1 cP for the gas, oil, and water phases, respectively. The phase relative permeabilities are represented using Brooks Corey model as follows,
\begin{equation}
\begin{aligned}
k_{rg}(S_{g}) = k_{rg,max}\left( \frac{S_{g}-S_{gr}}{1-Sgr-Sor-Swr}\right)^{n_{g}},\\
k_{ro}(S_{o}) = k_{ro,max}\left( \frac{S_{o}-S_{or}}{1-Sgr-Sor-Swr}\right)^{n_{o}},\\
k_{rw}(S_{w}) = k_{rw,max}\left( \frac{S_{w}-S_{wr}}{1-Sgr-Sor-Swr}\right)^{n_{w}},
\label{eqn:relperm}
\end{aligned}
\end{equation}
with endpoints $k_{rg,max}$ = 0.6, $k_{ro,max}$ = 0.7, $k_{rw,max}$ = 0.8, $S_{gr}$ = 0.1, $S_{or}$ = 0.15, $S_{wr}$ = 0.2, and model exponents $n_{g}$ = 1.5, $n_{o}$ = 1.2, $n_{w}$ = 2.0. Additionally, the gas-oil ($P_{cgo}$) and oil-water ($P_{cow}$) capillary pressures are also assumed to follow the Brook-Corey relationship,
\begin{equation}
\begin{aligned}
P_{cgo} = P_{tg}\left(\frac{1-S_{or}}{S_{o}-S_{or}} \right)^{\lambda_{go}},\\
P_{cow} = P_{to}\left(\frac{1-S_{wr}}{S_{w}-S_{wr}} \right)^{\lambda_{ow}}.
\end{aligned}
\end{equation}
Here, the entry pressures for the gas and water phases are $P_{tg}$ = 5 psi and $P_{to}$ = 10 psi, respectively. The model exponents are taken to be $\lambda_{go}$ = 0.5 and $\lambda_{ow}$ = 0.25. Figure \ref{fig:relcap} shows the phase relative permeabilities and capillary pressures obtained using the aforementioned Brooks-Corey model parameters.
The solution gas $R_{g}$ is assumed to be a monotonic function of pressure described by,
\begin{equation}
R_{g}(P_{o}) = 1 - exp^{-\beta (P_{o}-P_{d})}.
\end{equation}
Here, $\beta$ = 5$\times$10$^{-4}$ and $P_{d}$ = 1000 psi is the dew point pressure. Figure \ref{fig:solgas} shows the variation of solution gas as a function of pressure.
The initial conditions for the pressure and saturations are taken to be $P_{o}^{0}$ = 2500 psi, $S_{g}^{0}$ = 0.2, and $S_{o}^{0}$ = 0.55. We employ a pressure and water flux ($u_{w} \cdot n$) specified boundary conditions at the top right and bottom left corners of the computational domain to mimic pressure specified production at 2500 psi and water rate specified injection wells 1 STB/day. A no flow condition is assumed on the rest of the reservoir domain boundary. In what follows, we present a comparison of the oil, gas and water production rates and cumulative recoveries against a benchmark fine scale problem for both adaptive homogenization and multiscale basis approaches.

\subsection{Adaptive Numerical Homogenization}
\begin{figure}[H]
\begin{center}
\includegraphics[width=7.5cm,trim=5cm 10cm 4cm 9cm, clip]{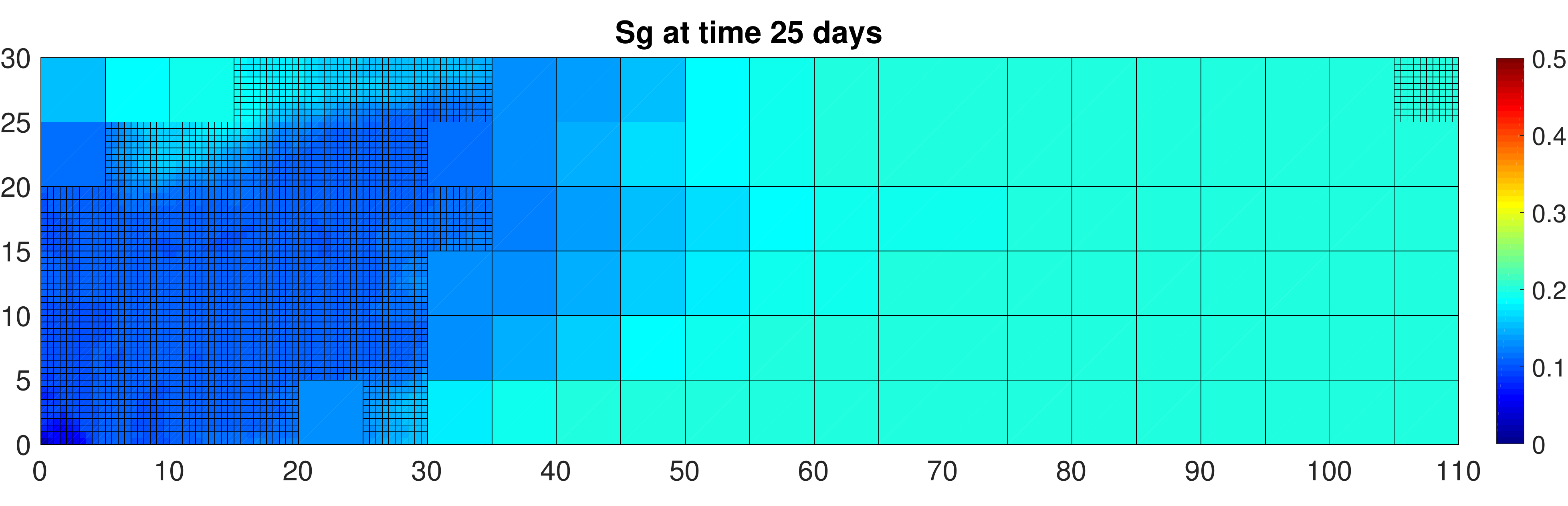}
\includegraphics[width=7.5cm,trim=5cm 10cm 4cm 9cm, clip]{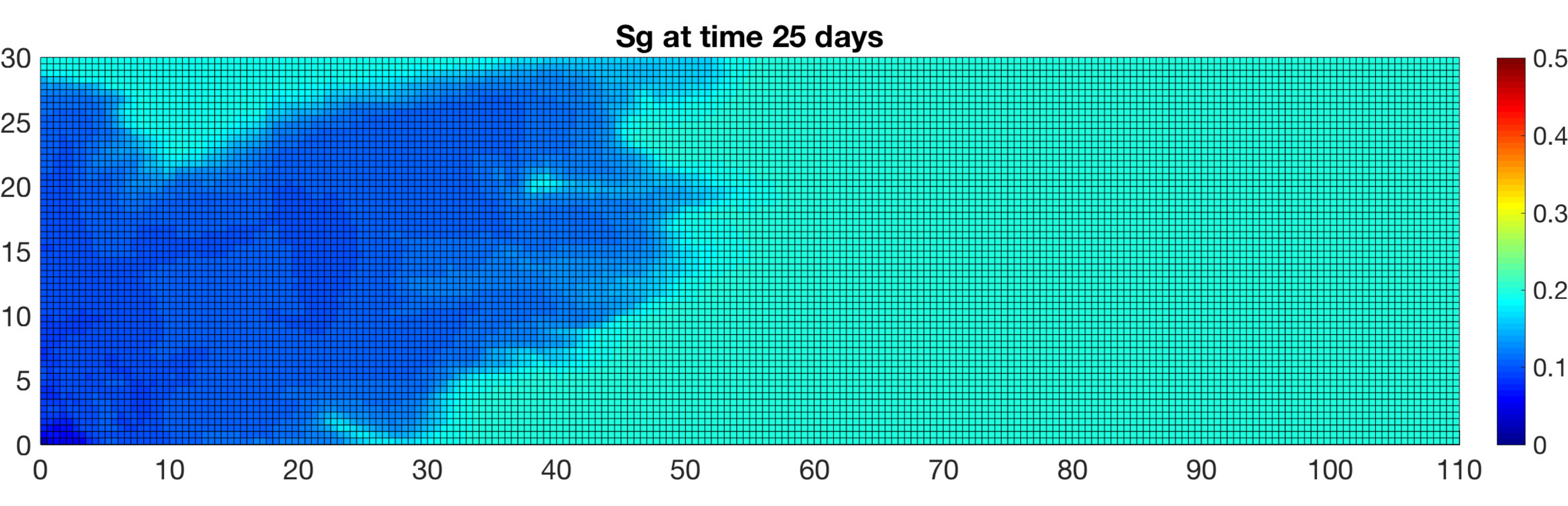}\\
\includegraphics[width=7.5cm,trim=5cm 10cm 4cm 9cm, clip]{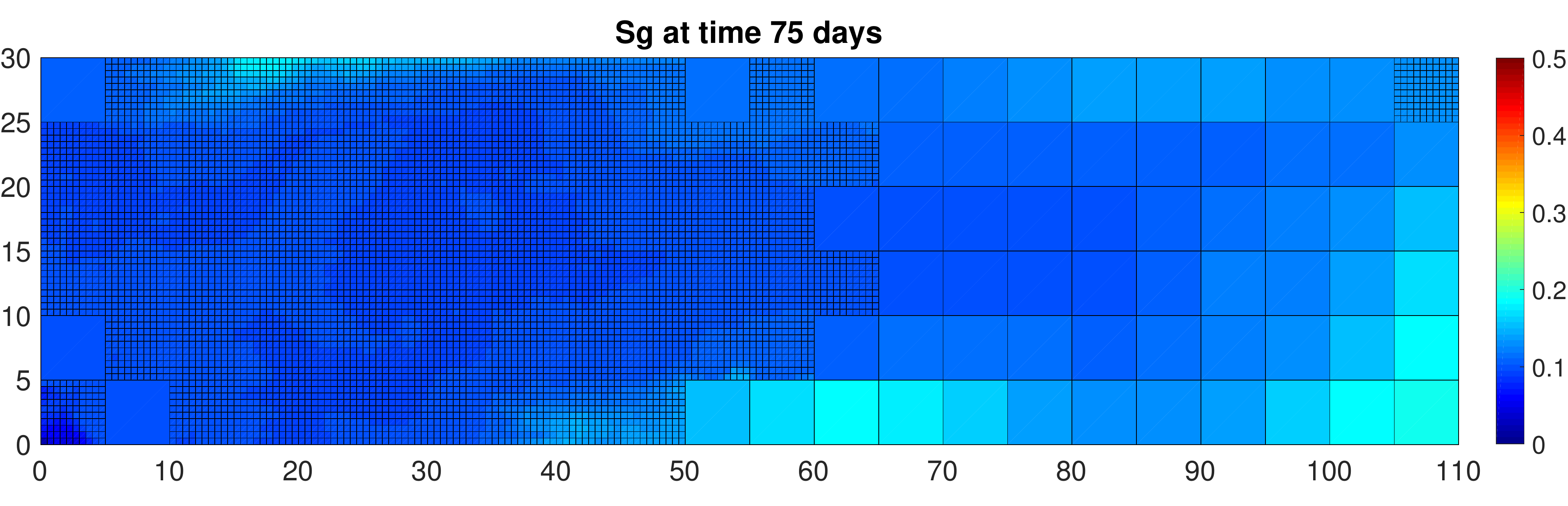}
\includegraphics[width=7.5cm,trim=5cm 10cm 4cm 9cm, clip]{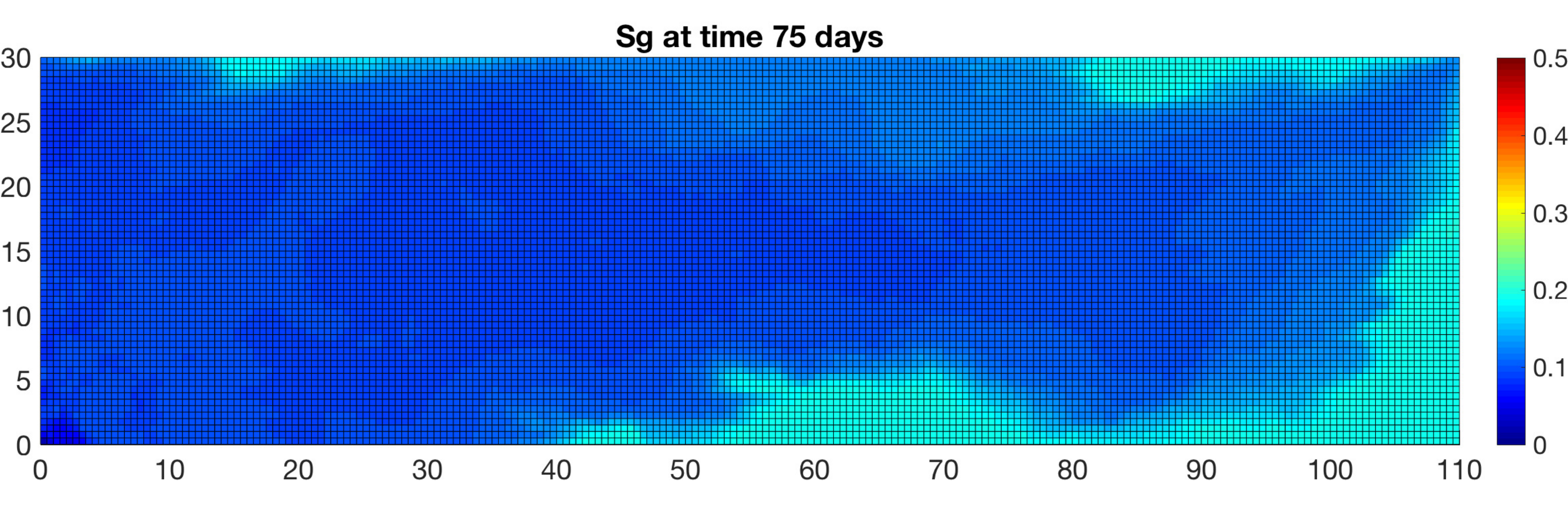}\\
\caption{Adaptive homogenization (left) and fine (right) scale gas saturation profiles}
\label{fig:Sg}
\end{center}
\end{figure}

\begin{figure}[H]
\begin{center}
\includegraphics[width=7.5cm,trim=5cm 10cm 4cm 9cm, clip]{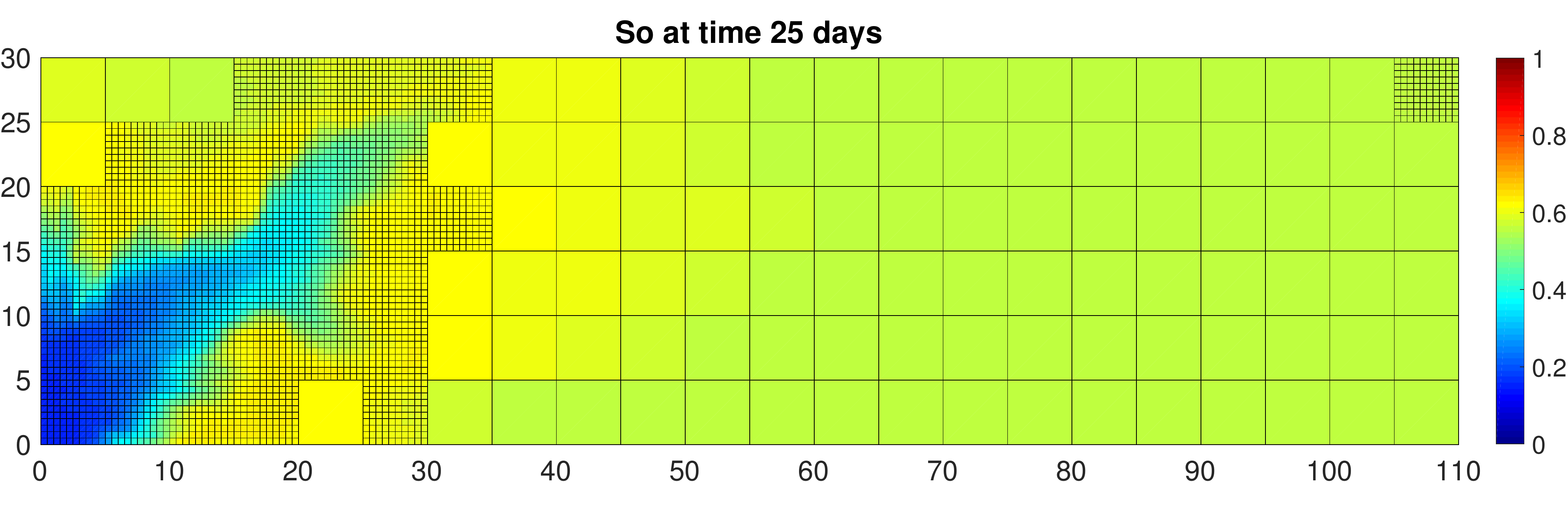}
\includegraphics[width=7.5cm,trim=5cm 10cm 4cm 9cm, clip]{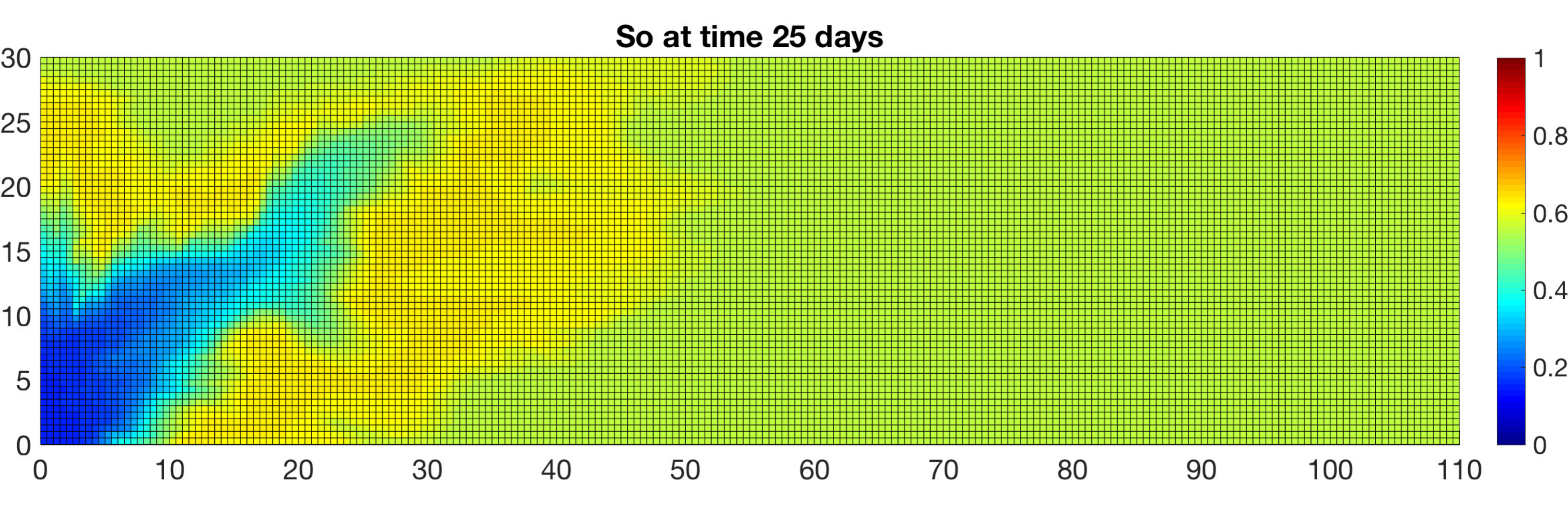}\\
\includegraphics[width=7.5cm,trim=5cm 10cm 4cm 9cm, clip]{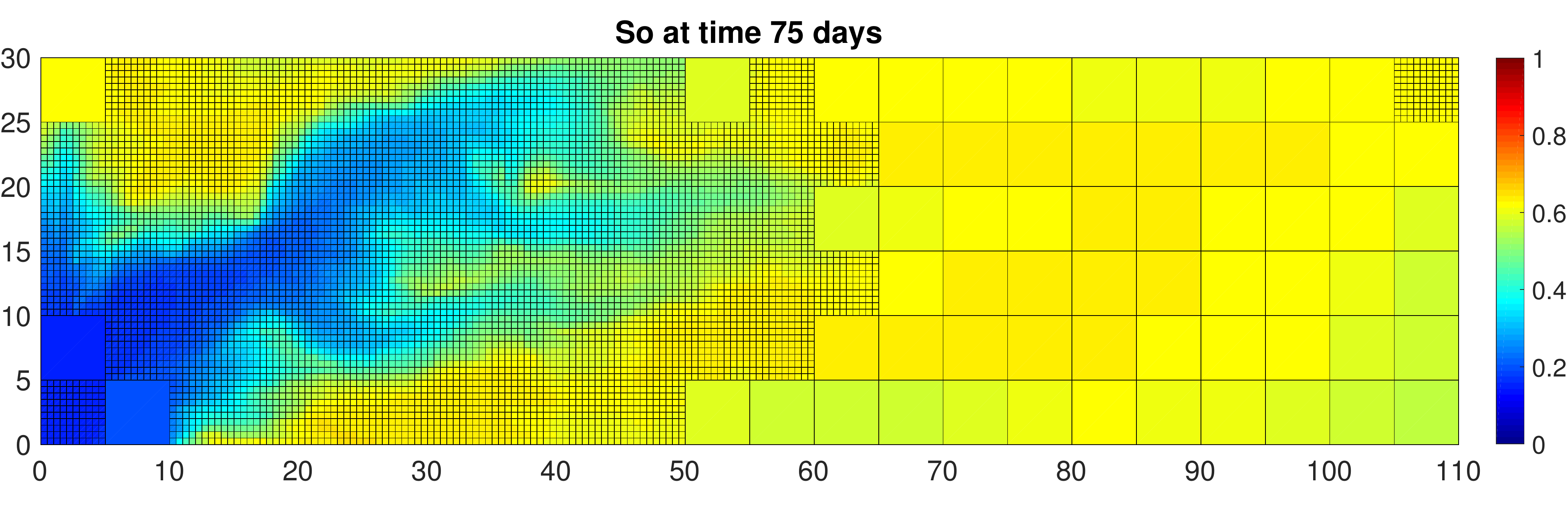}
\includegraphics[width=7.5cm,trim=5cm 10cm 4cm 9cm, clip]{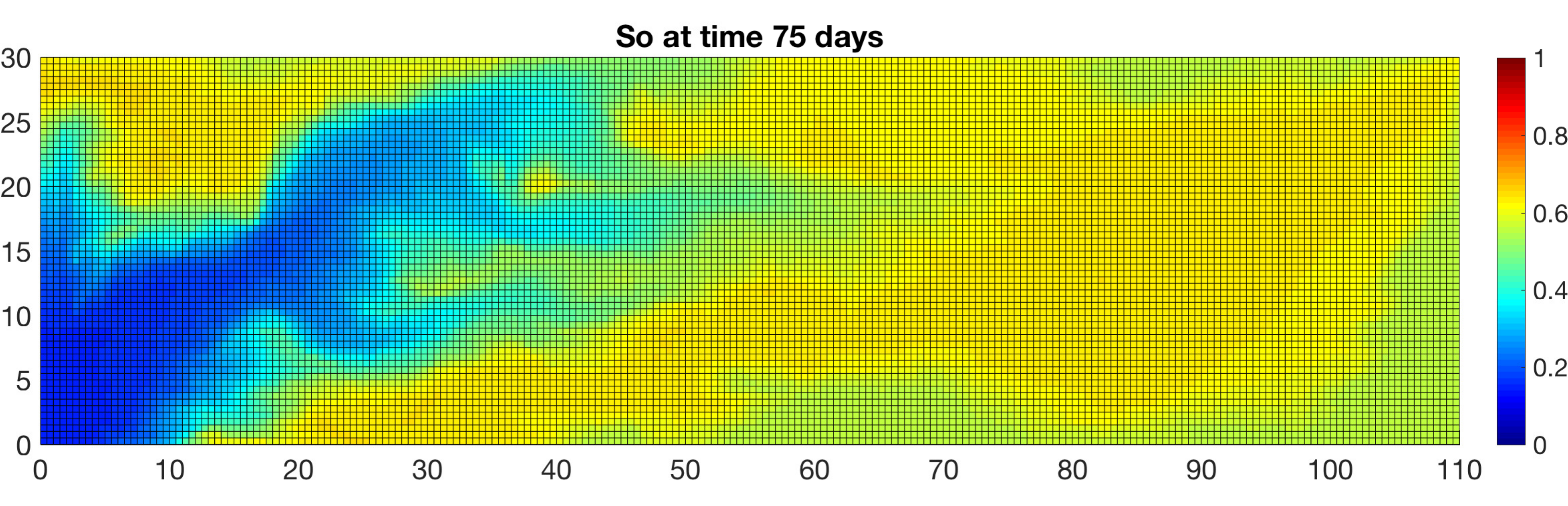}\\
\caption{Adaptive homogenization (left) and fine (right) scale oil saturation profiles}
\label{fig:So}
\end{center}
\end{figure}

\begin{figure}[H]
\begin{center}
\includegraphics[width=7.5cm,trim=5cm 10cm 4cm 9cm, clip]{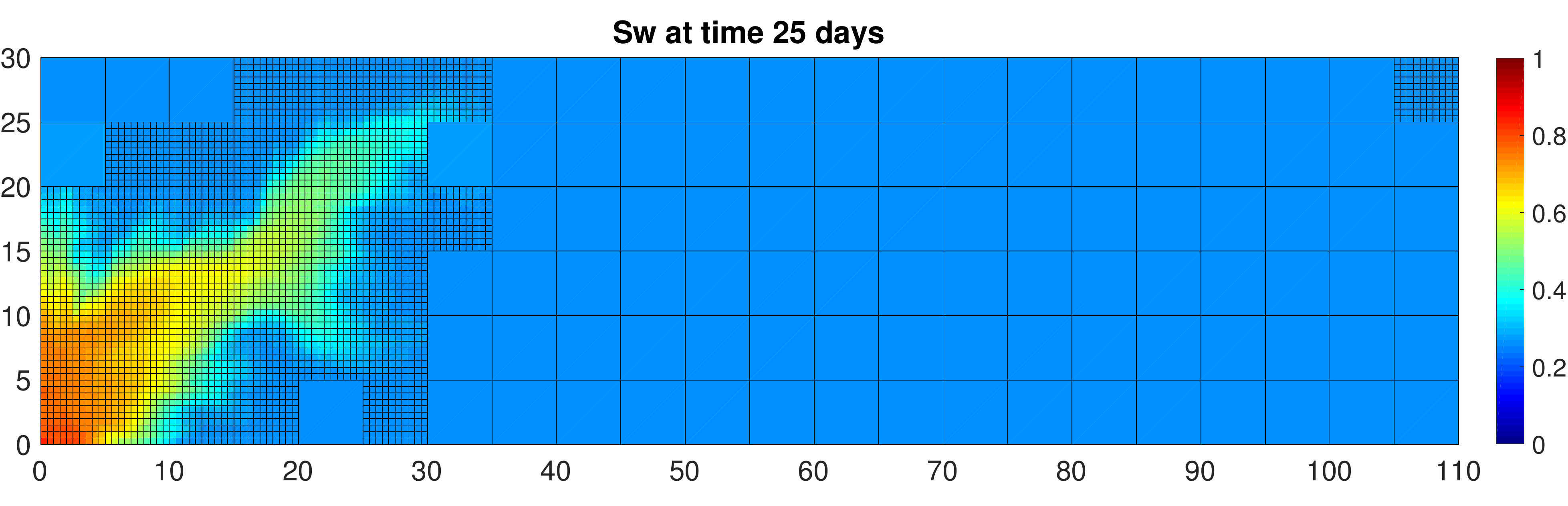}
\includegraphics[width=7.5cm,trim=5cm 10cm 4cm 9cm, clip]{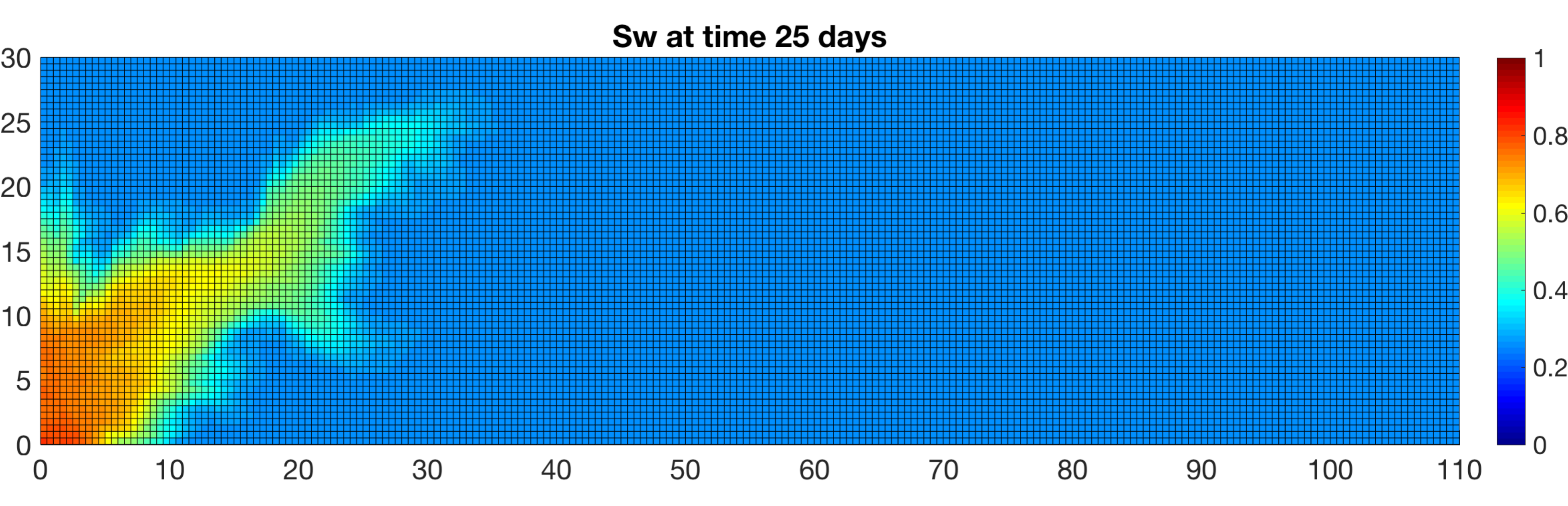}\\
\includegraphics[width=7.5cm,trim=5cm 10cm 4cm 9cm, clip]{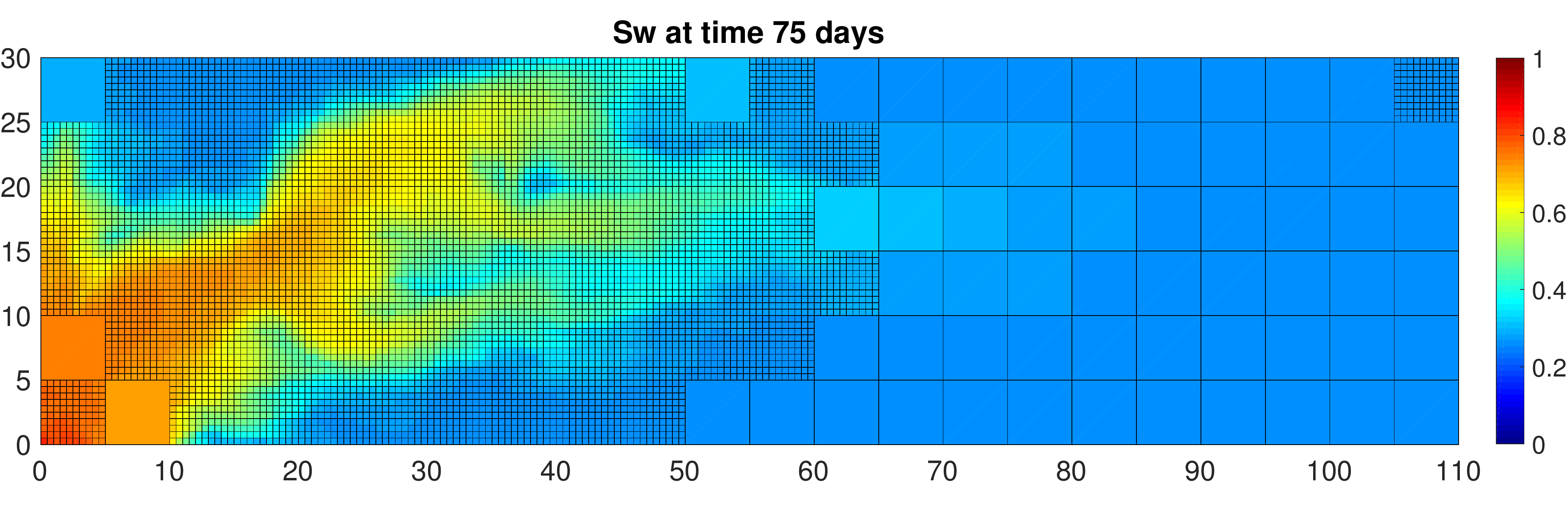}
\includegraphics[width=7.5cm,trim=5cm 10cm 4cm 9cm, clip]{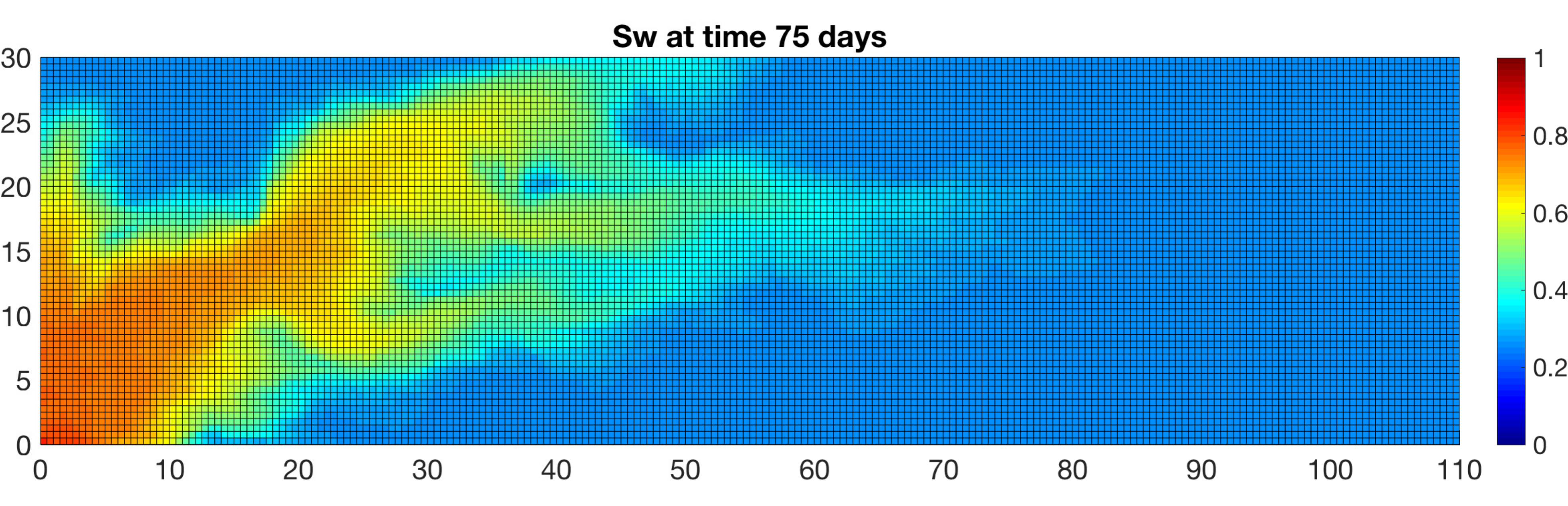}\\
\caption{Adaptive homogenization (left) and fine (right) scale water saturation profiles}
\label{fig:Sw}
\end{center}
\end{figure}

\begin{figure}[H]
\begin{center}
\includegraphics[width=7.5cm,trim=1.5cm 1cm 1cm 0cm, clip]{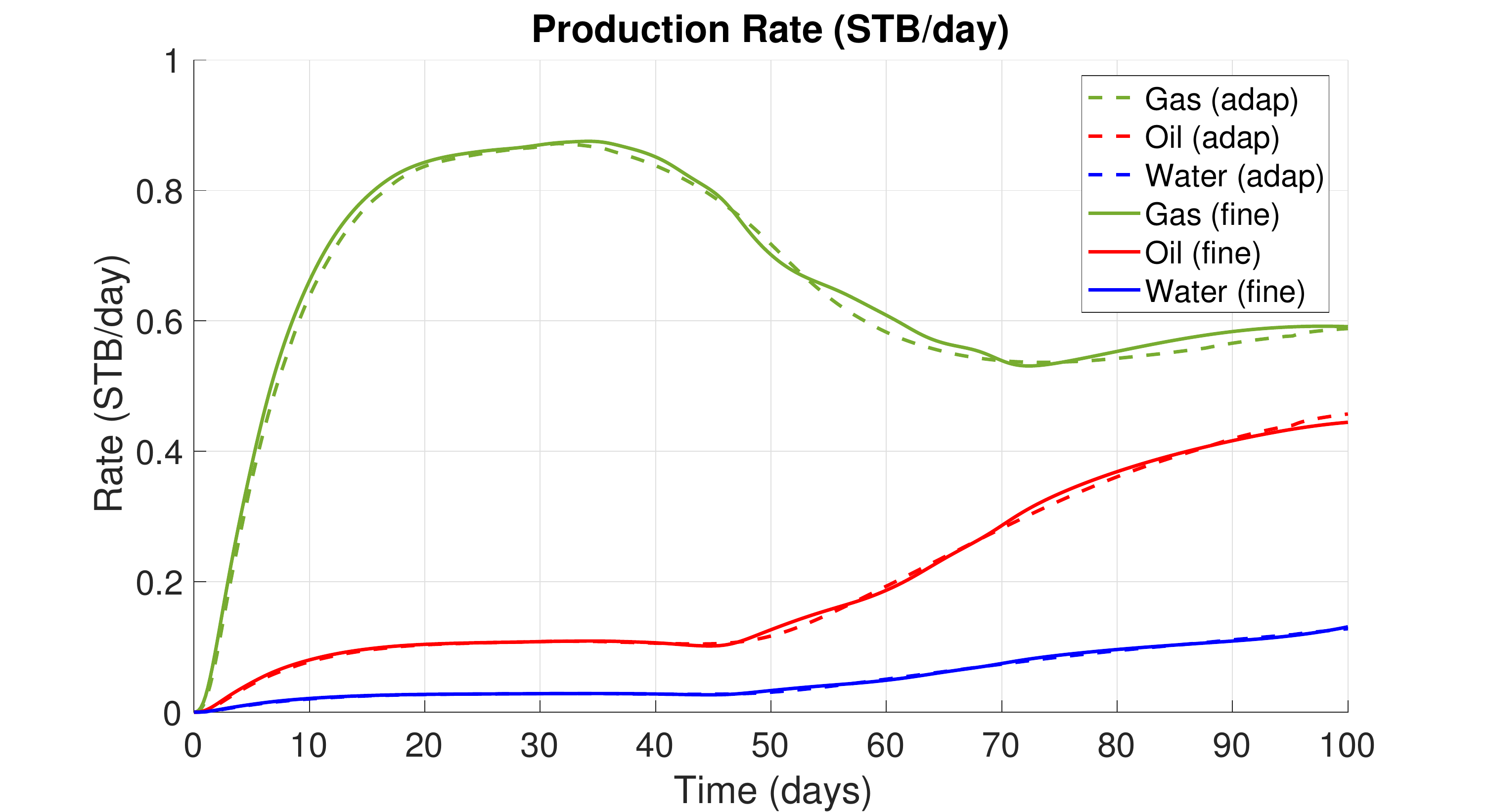}
\includegraphics[width=7.5cm,trim=1.5cm 1cm 1cm 0cm, clip]{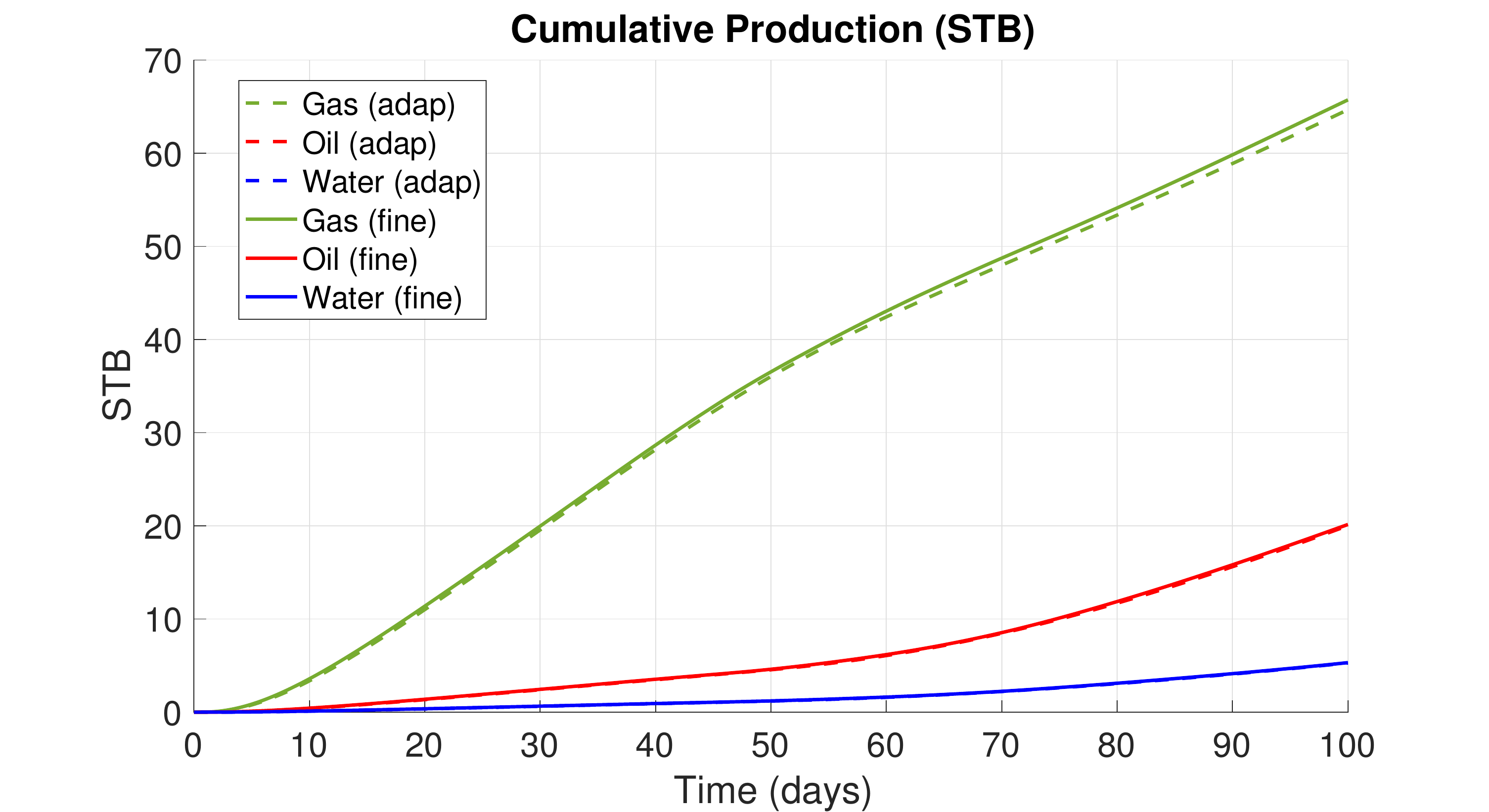}
\caption{Gas, oil, and water rates and cumulative recoveries at the production well for fine scale and adaptive homogenization.}
\label{fig:recov}
\end{center}
\end{figure}

In this subsection we present a comparison between a benchmark fine scale simulation and adaptive numerical homogenization approach. The tolerance for the adaptivity criterion for dynamic mesh refinement, described in subsection \ref{subsec:criterion}, is $\epsilon_{adap} =$ 0.05. Figures \ref{fig:Sg}, \ref{fig:So}, and \ref{fig:Sw} compare the gas, oil, and water saturation profiles, respectively between the adaptive and fine scale simulations at time 25 and 75 days. The adaptivity criterion identifies the water saturation front as the transient region and adaptively refines this region of interest while using coarse scale effective properties away from the transient region.
We note that the gas saturation front , in Figure \ref{fig:Sg} moves faster than the oil and water saturation fronts, Figures \ref{fig:So} and \ref{fig:Sw}, owing its high phase mobility. However, identifying the water saturation front as the transient region allows us to capture the flow physics without requiring local mesh refinement (or enrichment) on a larger computational domain. Furthermore, we are also able to accurately capture the oil bank, region of high oil saturation away from the water saturation front in Figure \ref{fig:So}, formed due to large, mobile, initial oil saturations.
Figure \ref{fig:recov} shows a comparison between adaptive and fine scale simulations for the gas, oil, and water rates and cumulative recoveries at the production well.
The results show that the adaptive homogenization approach is in good agreement with the benchmark fine scale simulations for the given reservoir properties and model parameters.

\subsection{Generalized mixed multiscale method}
\begin{figure}[H]
\begin{center}
\includegraphics[width=75mm,trim=1.5cm 1cm 1cm 0cm,scale=0.5]{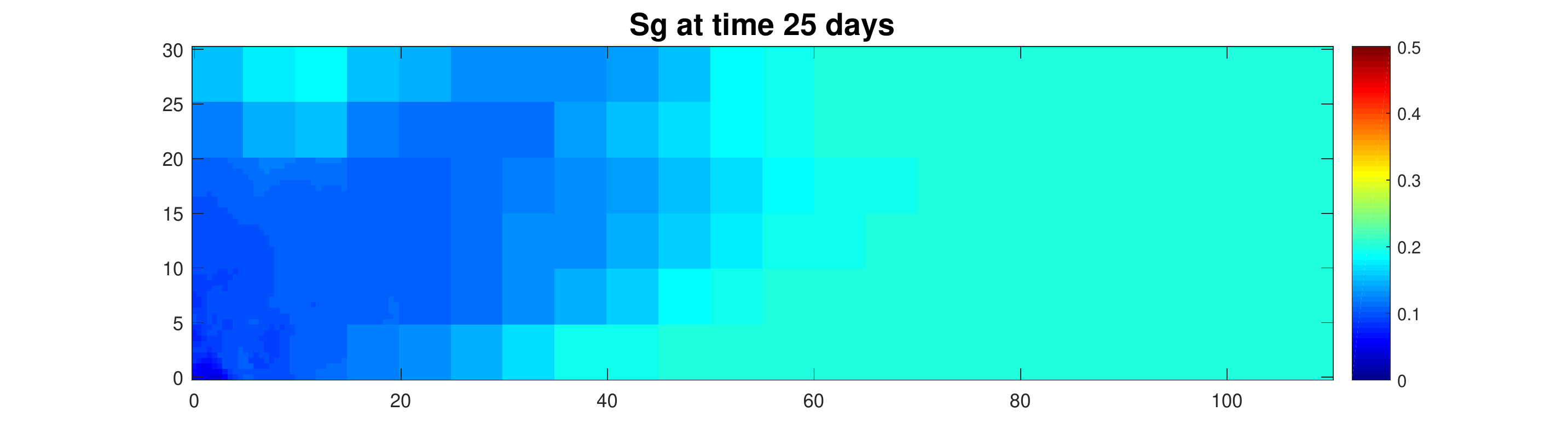}
\includegraphics[width=75mm,trim=1.5cm 1cm 1cm 0cm,scale=0.5]{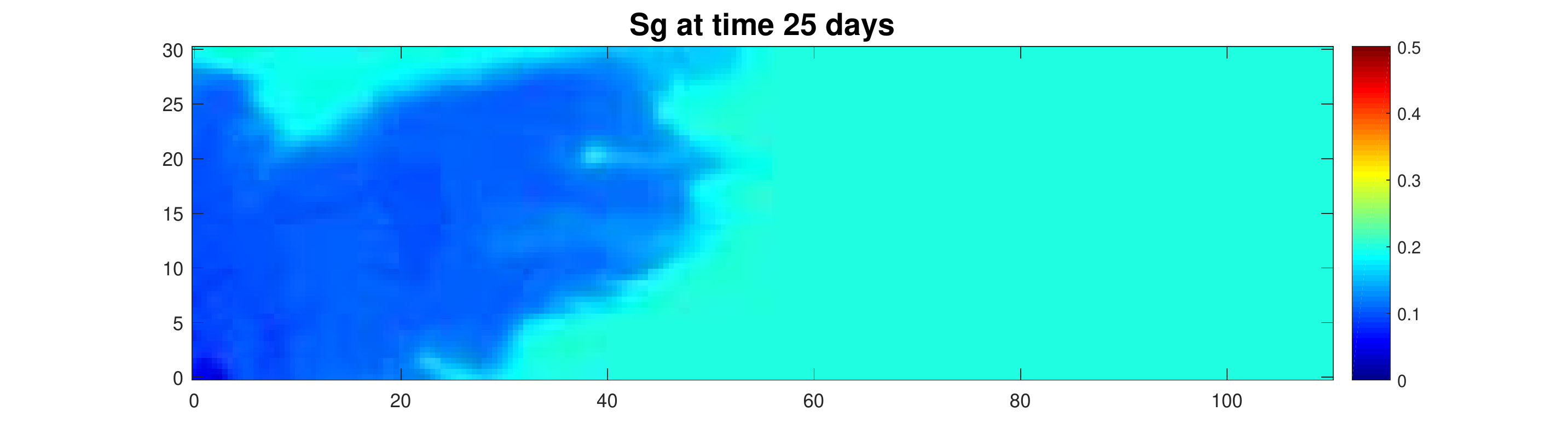}\\
\includegraphics[width=75mm,trim=1.5cm 1cm 1cm 0cm,scale=0.5]{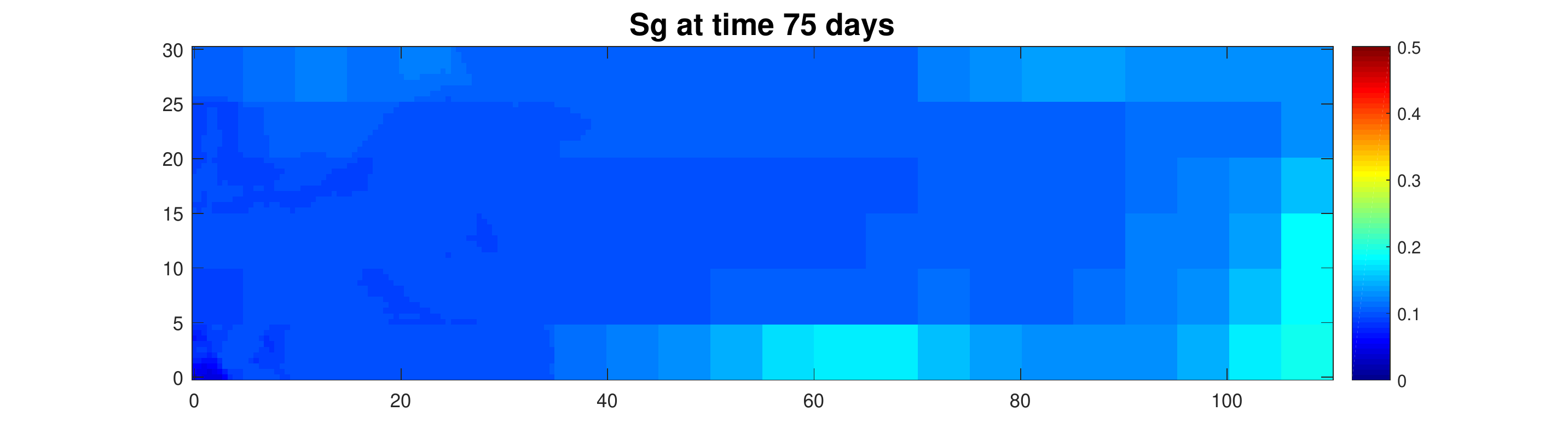}
\includegraphics[width=75mm,trim=1.5cm 1cm 1cm 0cm,scale=0.5]{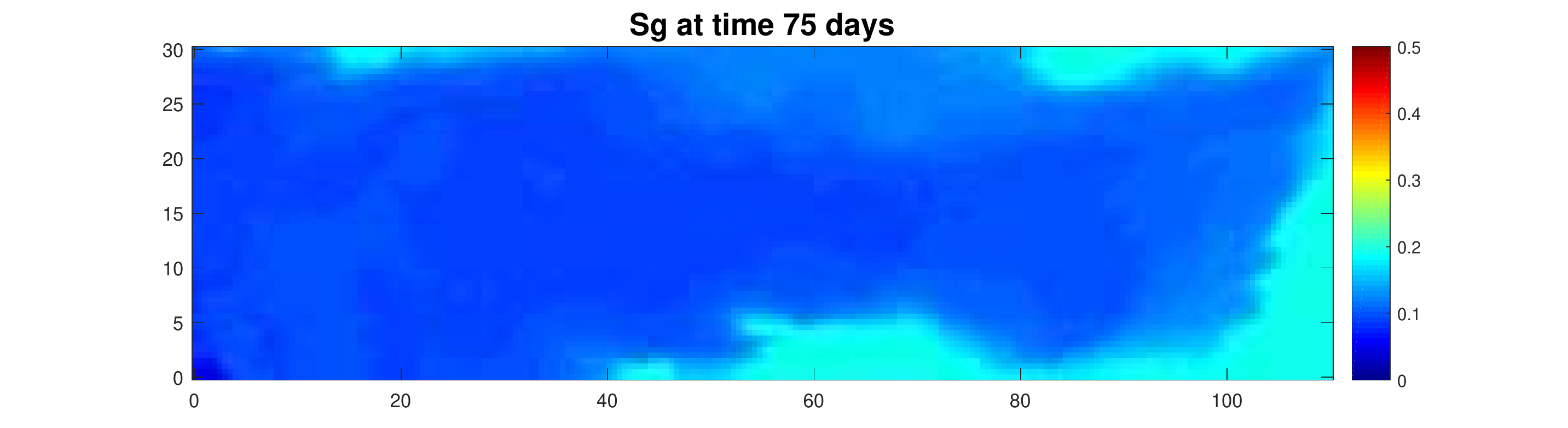}\\
\caption{Generalized multiscale (left) and fine (right) scale gas saturation profiles}
\label{fig:Sg_ms}
\end{center}
\end{figure}

\begin{figure}[H]
\begin{center}
\includegraphics[width=75mm,trim=1.5cm 1cm 1cm 0cm,scale=0.5]{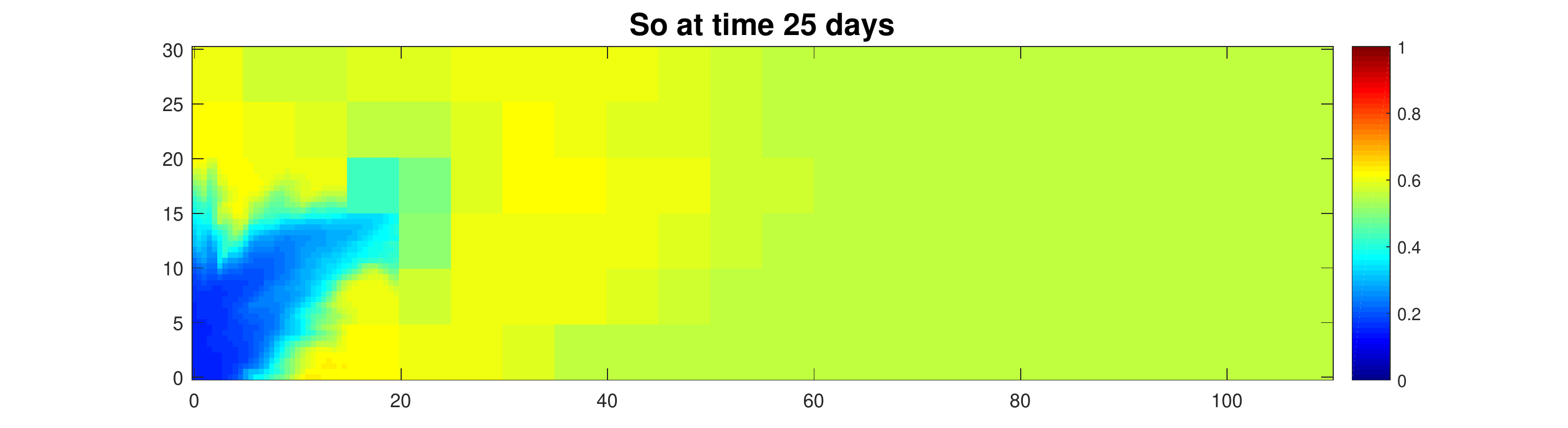}
\includegraphics[width=75mm,trim=1.5cm 1cm 1cm 0cm,scale=0.5]{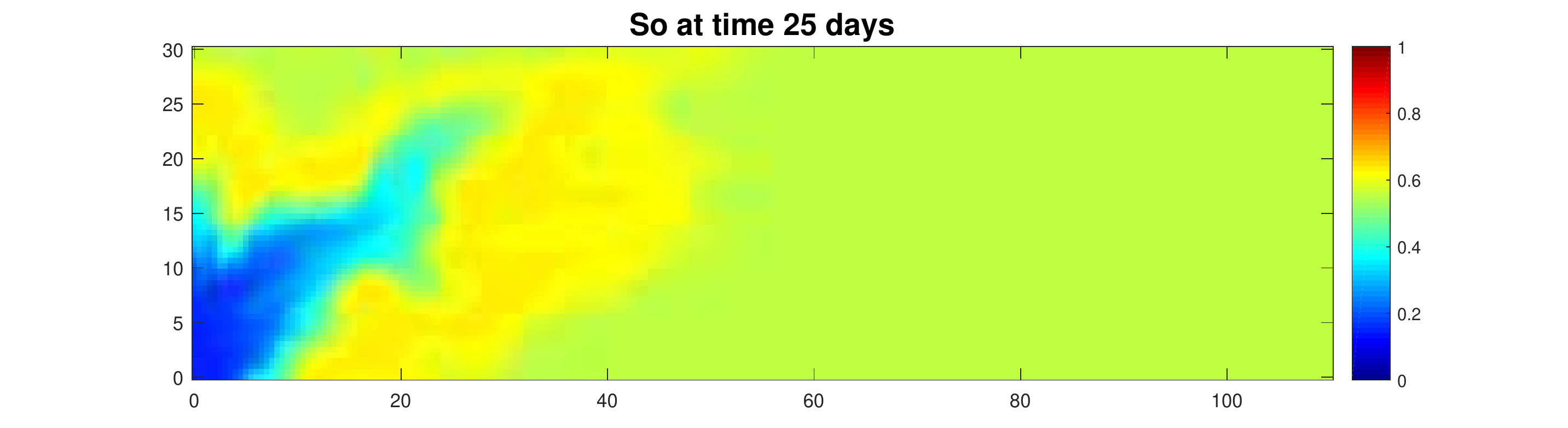}\\
\includegraphics[width=75mm,trim=1.5cm 1cm 1cm 0cm,scale=0.5]{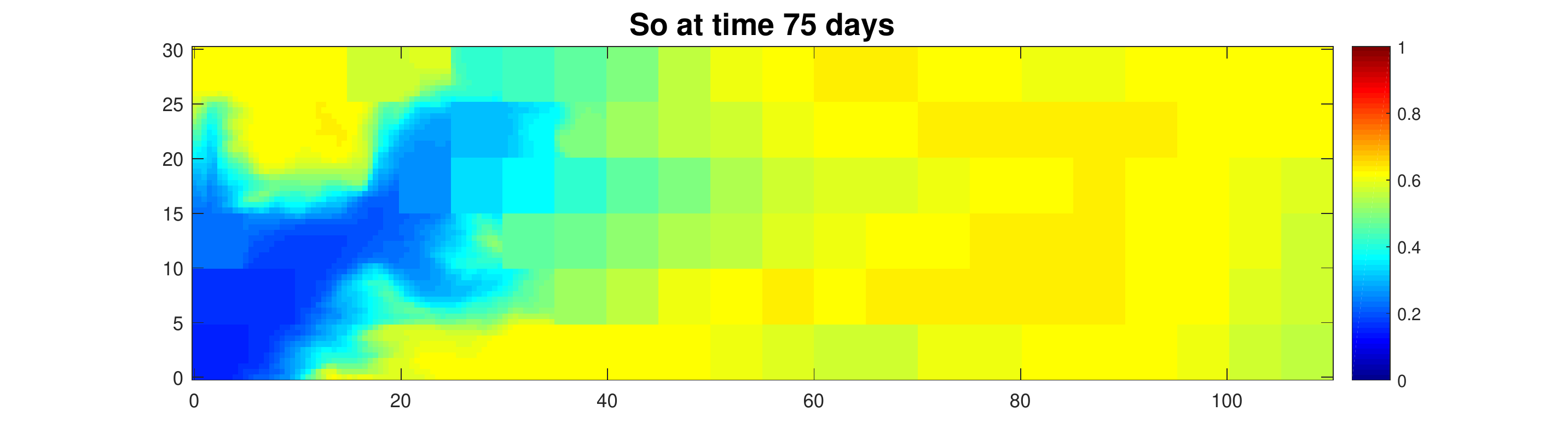}
\includegraphics[width=75mm,trim=1.5cm 1cm 1cm 0cm,scale=0.5]{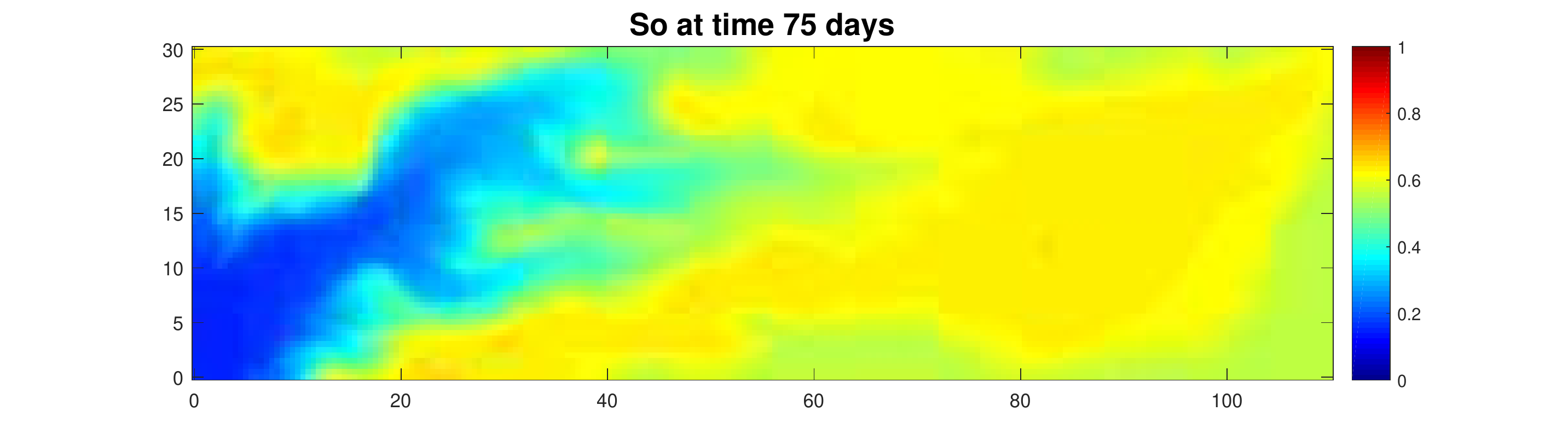}\\
\caption{Generalized multiscale (left) and fine (right) scale oil saturation profiles}
\label{fig:So_ms}
\end{center}
\end{figure}

\begin{figure}[H]
\begin{center}
\includegraphics[width=75mm,trim=1.5cm 1cm 1cm 0cm,,scale=0.5]{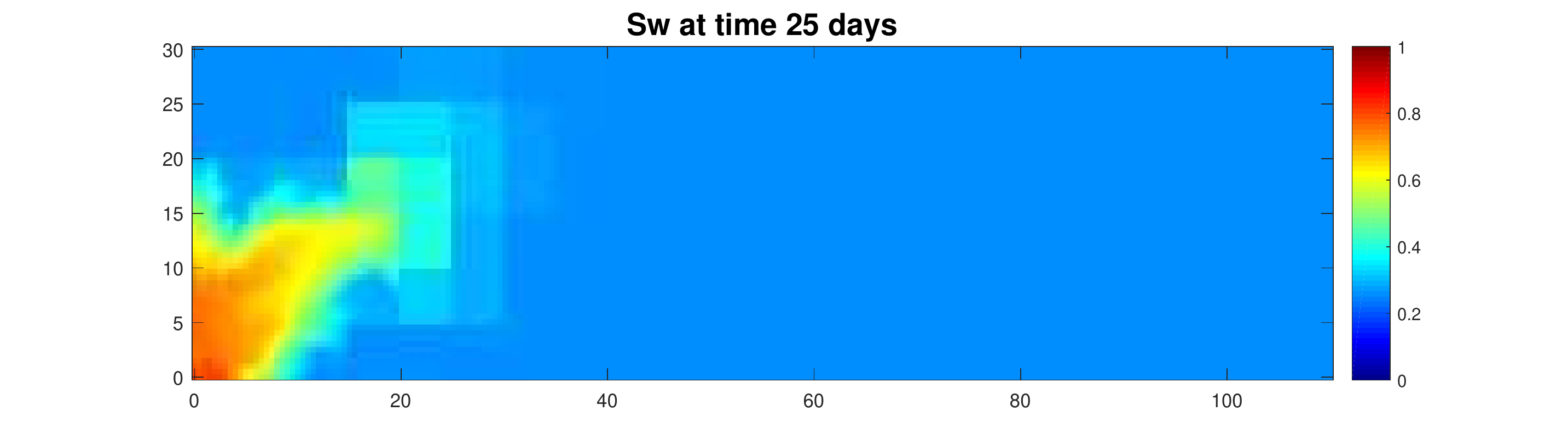}
\includegraphics[width=75mm,trim=1.5cm 1cm 1cm 0cm,,scale=0.5]{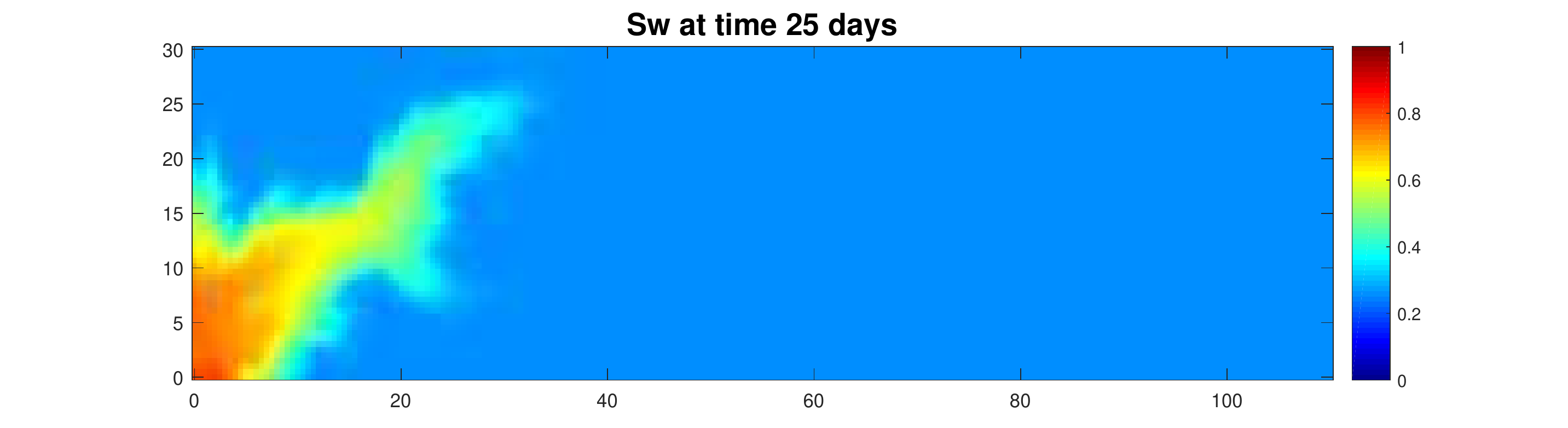}\\
\includegraphics[width=75mm,trim=1.5cm 1cm 1cm 0cm,,scale=0.5]{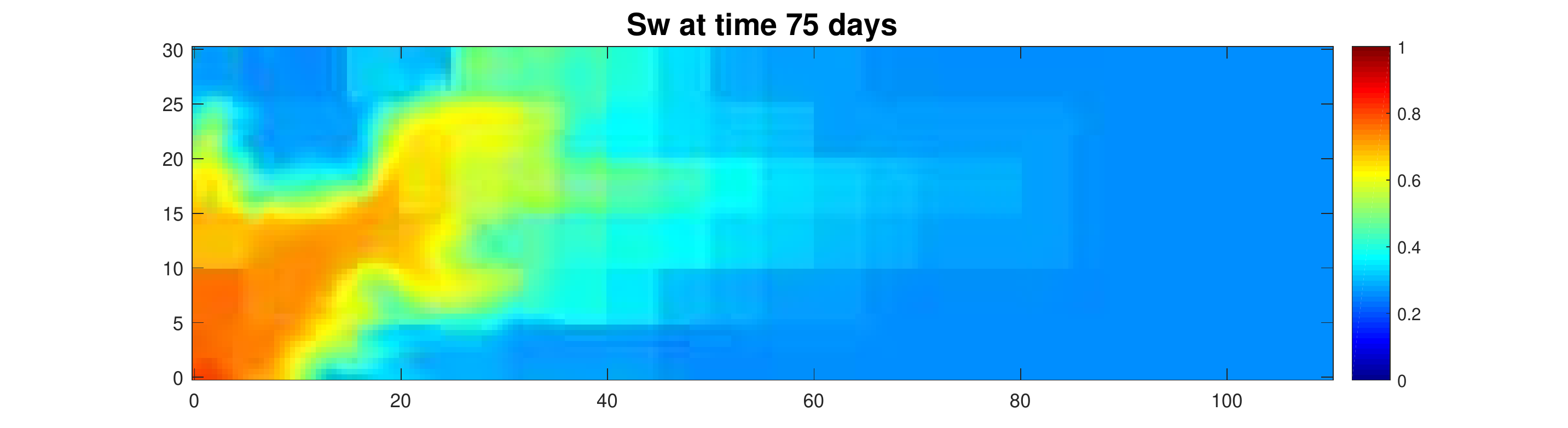}
\includegraphics[width=75mm,trim=1.5cm 1cm 1cm 0cm,,scale=0.5]{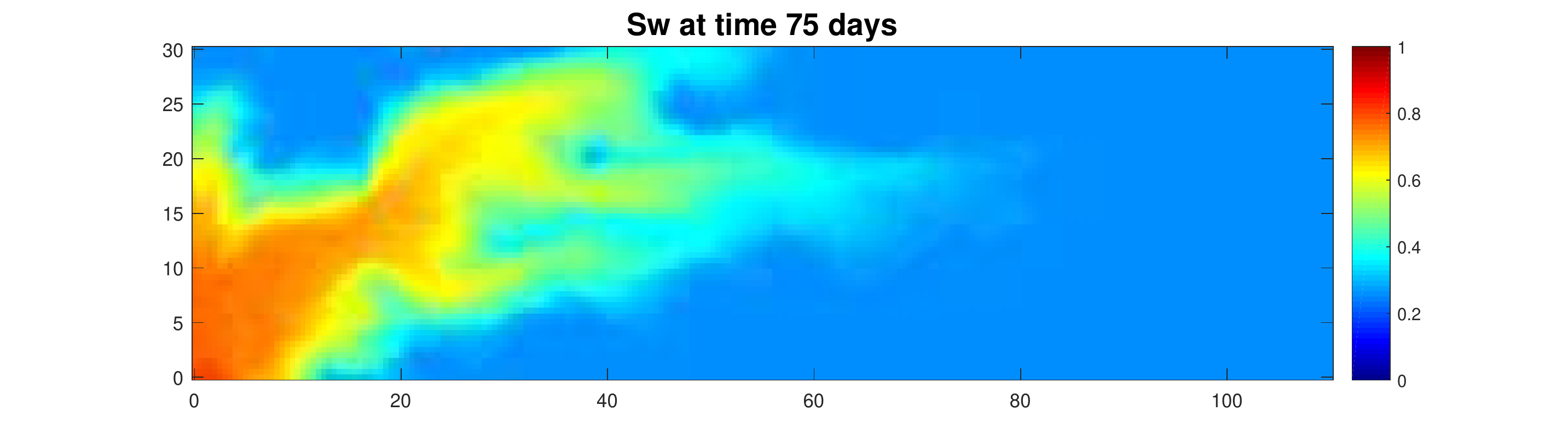}\\
\caption{Generalized multiscale (left) and fine (right) scale water saturation profiles}
\label{fig:Sw_ms}
\end{center}
\end{figure}

\begin{figure}[H]
\begin{center}
\includegraphics[width=7.5cm,trim=1.5cm 1cm 1cm 0cm, clip]{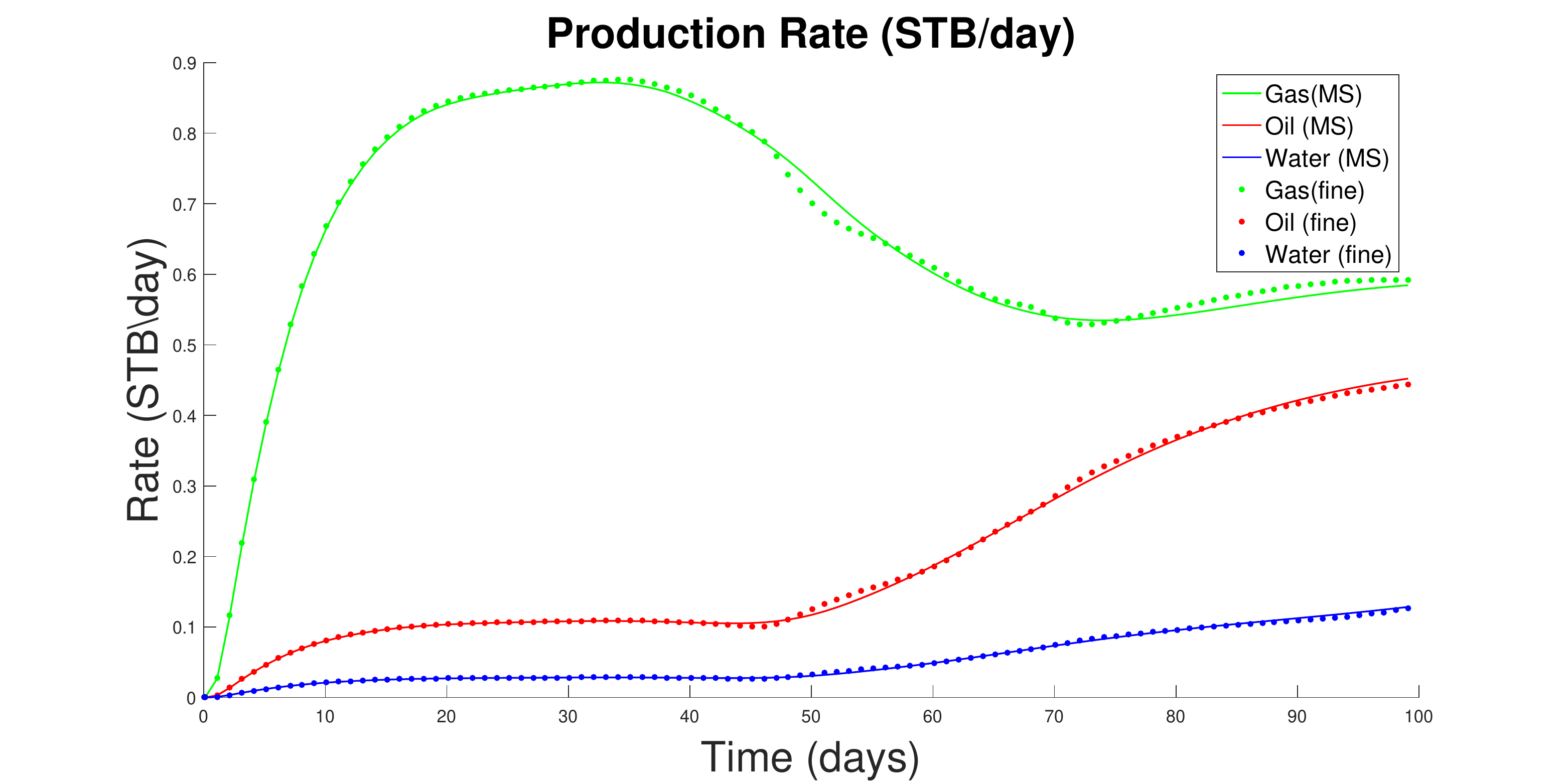}
\includegraphics[width=7.5cm,trim=1.5cm 1cm 1cm 0cm, clip]{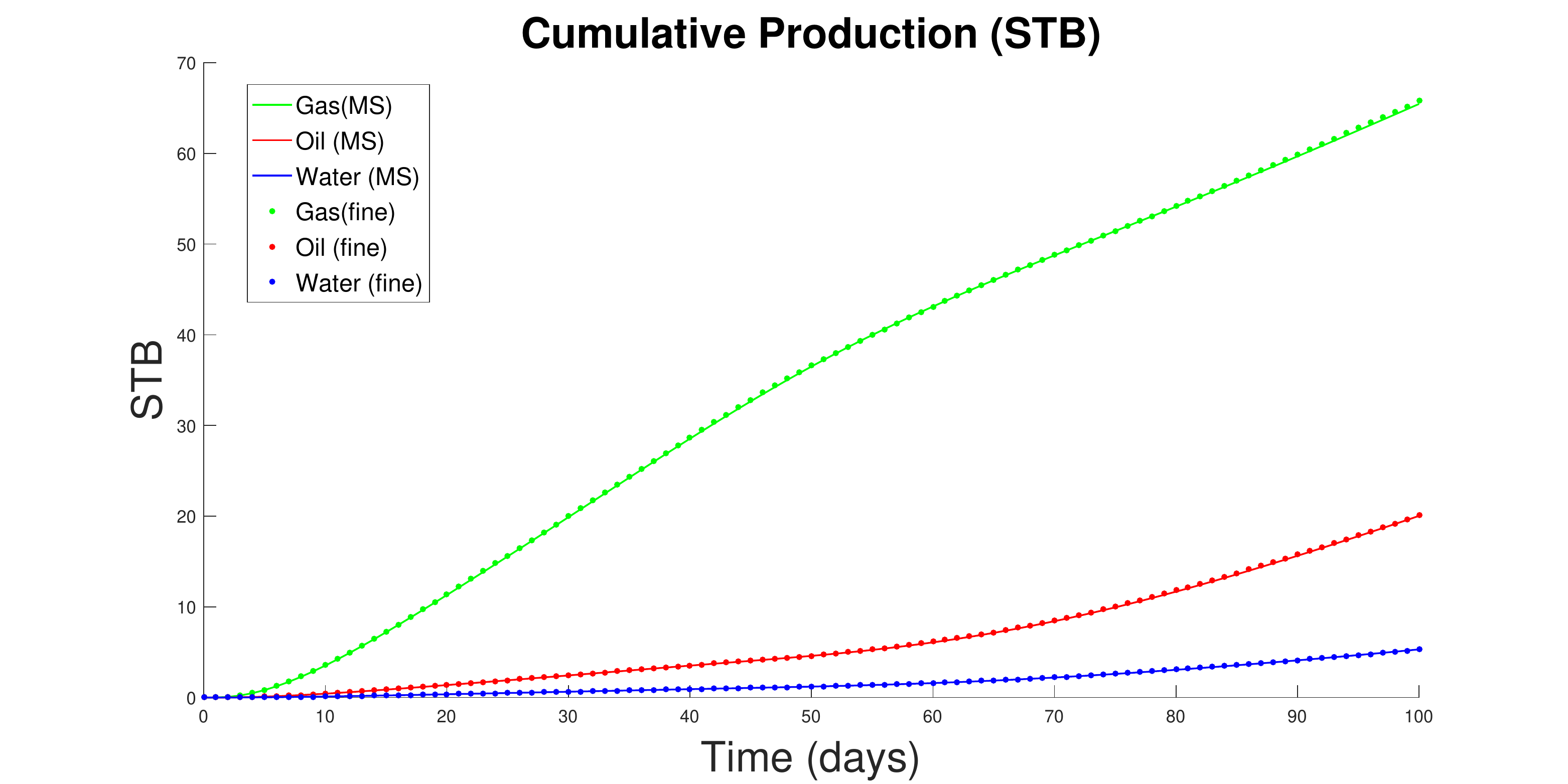}
\caption{Gas, oil, and water rates and cumulative recoveries at the production well for fine scale and generalized multiscale.}
\label{fig:recov_ms}
\end{center}
\end{figure}

In this subsection, we will present a numerical result of the generalized mixed multiscale method for black oil equation. We will use the adaptivity criterion described in \ref{sec:adaptive_update} with $\theta=.04$.  In Figure \ref{fig:Sw_ms}, \ref{fig:So_ms} , \ref{fig:Sg_ms}, we will show the snapshots of the gas, oil and water saturation for both multiscale solution and the fine grid reference solution at time 25 days and 75 days. We can see that even we are using coarse grid in most of the regions, we still obtain a good coarse grid approximation of the reference solution.
We observe that the fine grid refinement is mostly concentrated at the front of the water saturation front since the mass transfer in the water phase is the dominating part in the total mass transfer. In most of the regions, we are still using coarse grid piesewise constant basis functions. The adaptive approach keeps the problem size small to increase the computational speed. To improve the accuracy, we can either enriching multiscale basis function to the coarse grid space, instead of using piecewise constant function only, or using smaller adaptive tolerance. Using enriching multiscale basis function gives a better global approximation for the pressure and velocity. Using smaller adaptive tolerance can give us a finer resolution for the saturations, the pressure and the velocity locally.
 In Figure \ref{fig:recov_ms}, we show show a comparison between adaptive and fine scale simulations for the gas, oil, and water rates and cumulative recoveries at the production well. We can notice that the production rate for the multiscale method is close to the reference production rate. The production curve is a little bit smoother that the reference production curve due to the coarse grid effect. Looking at the cumulative production, we find two curves are nearly the same.

\section{Summary} \label{sec:summary}

We presented two approaches for model order reduction namely: (1) the adaptive numerical homogenization and (2) the generalized multiscale basis functions. An expanded mixed variational formulation is used by both approaches to separate scale of the flow and transport problems. The adaptive numerical homogenization approach performs local numerical homogenization to pre-compute effective properties (permeability, porosity, relative permeability, and capillary pressure) and dynamic mesh refinement to fine scale during simulation runtime based upon an adaptivity criterion. On the other hand, the generalized multiscale basis functions approach construct a multiscale velocity basis functions as an offline (or pre-compute) step. We describe a benchmark fine scale, black oil model problem to compare the accuracy of each of these two approaches. The numerical results indicate that both approaches accurately capture the multiphase flow physics while reducing the problem size compared to the benchmark fine scale simulation. This benchmark study enables us to rigorously assess, in the near future, the computational efficiencies required of each of the aforementioned model order reduction techniques.

\bibliographystyle{plain}  
\bibliography{ref}

\end{document}